\documentclass[12pt,a4paper]{article}
\usepackage[includeheadfoot,left=2.5cm,right=2.5cm,top=2.5cm,bottom=2.5cm,head=15pt,
bindingoffset=0pt]{geometry}
\usepackage{preamble}

\newcommand{\Ranklocus}[4]{\varUpsilon_{#1}^{#2,#3}(#4)}

\newcommand{\cactus}[2]{\mathfrak{K}_{#1}\left( #2 \right)}

\newcommand{\annihp}{\ccA{}nn(p)}
\newcommand{\indexa}{i}
\newcommand{\indexb}{j}
\newcommand{\indexc}{k}

\title{Secant varieties to high degree Veronese reembeddings, catalecticant matrices and smoothable Gorenstein schemes}
\author{\nisiabu \and \JaBu\thanks{Supported by Marie Curie Outgoing Fellowship ``Contact Manifolds''.}}

\date{22 November 2011}

\begin{document}

\maketitle

\begin{abstract}
We study the secant varieties of the Veronese varieties and of Veronese reembeddings of a smooth projective variety.
We give some conditions, under which these secant varieties are set-theoretically cut out by determinantal equations.
More precisely, they are given by minors of a catalecticant matrix.
These conditions include the case when the dimension
of the projective variety is at most 3 and the degree of reembedding is sufficiently high.
This gives a positive answer to a set-theoretic version of a question of Eisenbud in dimension at most 3.
For dimension four and higher
we produce plenty of examples when the catalecticant minors are not enough
  to set-theoretically define the secant varieties to high degree Veronese varieties.
This is done by relating the problem to smoothability of certain zero-dimensional Gorenstein schemes.
\end{abstract}

\medskip
{\footnotesize
\noindent\textbf{e-mail addresses:} \verb|wkrych@mimuw.edu.pl|  and  \verb|jabu@mimuw.edu.pl|

\noindent\textbf{address:} Institut Fourier, 100 rue des Maths, BP 74, 38402 St Martin d'H\`eres Cedex, France

\noindent\textbf{keywords:}
secant variety, catalecticant minors, Veronese variety,
Veronese reembeddings, cactus variety, smoothable zero-dimensional schemes.

}

\section{Introduction}\label{section_intro}

Throughout the paper we work over the base field of complex numbers $\CC$.
We investigate the secant varieties to Veronese embeddings of projective space.
Despite this being a topic of a very intensive research
  (see \cite{iarrobino_kanev_book_Gorenstein_algebras},
       \cite{kanev_chordal_varieties},
       \cite{landsberg_ottaviani_equ_for_secants_to_Veronese_varieties},
       \cite{raicu_3_times_3_minors}
     and references therein),
  the defining equations of such varieties
  are hardly known except in few cases
  (see \cite[Table on p.2]{landsberg_ottaviani_equ_for_secants_to_Veronese_varieties}
    for the list of known cases and also the more recent \cite{raicu_3_times_3_minors},
    where the results for the second secant variety are improved).
There are however some equations known: these are the determinantal equations arising from catalecticant matrices.
See \cite{geramita_inverse_systems_of_fat_points}, \cite{geramita_catalecticant_varieties}
   for an overview on the ideals generated by catalecticant minors.
In general, these equations are not enough to define the secant variety, even set-theoretically,
but in the known examples the degree of the reembedding is relatively low.
In this paper we assume that the degree of the Veronese reembedding is sufficiently high.
With this assumption, we prove that there are few situations when the determinantal equations are sufficient to define
  the secant variety set-theoretically, see Theorem~\ref{thm_secants_are_cut_out_by_catalecticants}.
We also observe that even for high degrees the catalecticant minors are rarely sufficient
  to define the secant varieties, see Theorem~\ref{thm_counterexamples_to_Eisenbud}.

\subsection{Secant varieties to Veronese varieties and catalecticants}

Throughout the article $V$ and $W$ denote complex vector spaces, and $\PP V$ is the naive projectivisation of $V$.
For $v \in V \setminus \set{0}$,  by $[v]$ we mean the corresponding point in the projective space $\PP V$.
We state our first theorem and then we explain the notation used in the theorem in details.

\begin{thm}\label{thm_secants_are_cut_out_by_catalecticants}
  Let $r, n, d, \indexa$ be four integers and let $V\simeq \CC^{n+1}$ be a vector space.
  Let $\sigma_r(v_d (\PP V))$ be the $r$-th secant variety of $d$-th Veronese embedding of $\PP V \simeq \PP^n$.
  If $d \ge 2r$, $r \le \indexa \le d-r$
     and also either $r \le 10$ or $n \le 3$,
  then $\sigma_r(v_d (\PP V))$ is set-theoretically defined by $(r+1)\times(r+1)$ minors of the $\indexa$-th catalecticant matrix.
\end{thm}

The theorem extends to the case of $\sigma_r(v_d (X))$, where $X$ is a smooth projective variety of dimension $n$,
   see Corollary~\ref{cor_Eisenbuds_holds_for_X} below.
In Section~\ref{section_improving_bounds} we briefly discuss potential extensions of the range of the integers in the theorem.
We give the proof in Section~\ref{section_smoothability_of_Gorenstein}.

Throughout the article, for a subvariety $X \subset \PP W$, the \emph{$r$-th secant variety} $\sigma_r(X)$ is defined as
\[
  \sigma_r(X):= \overline{\bigcup_{\fromto{x_1}{x_r}\in X}  \langle \fromto{x_1}{x_r}\rangle}\subset \PP W
\]
where $\langle \fromto{x_1}{x_r}\rangle\subset \PP W$ denotes the linear span of the points $\fromto{x_1}{x_r}$
and the overline denotes Zariski closure.
In particular, $\sigma_1(X) = X$.

The $d$-th symmetric tensor power of $V$ is denoted $S^d V$.
The symmetric algebra of polynomial functions on $V$
 is denoted $\Sym V^* := \bigoplus_{d=0}^{\infty} S^d V^*$.

The \emph{$d$-th Veronese embedding} is denoted $v_d\colon \PP V \to \PP (S^d V)$
  and it maps $[v]$ for $v\in V$ to the class of the symmetric tensor $[v^d]$.

For $\indexa \in \setfromto{0,1}{d-1,d}$ we have a natural linear embedding
  \[
    \iota_{\indexa}\colon \PP (S^d V)  \to \PP(S^{\indexa} V \otimes S^{d- \indexa} V).
  \]
This comes from the embedding of the space of symmetric tensors into the space of partially symmetric tensors.
The map $\iota_{\indexa}$ is called the \emph{$\indexa$-flattening map}.
We define:
\begin{equation}\label{equ_define_ranklocus}
   \Ranklocus{r}{\indexa}{d-\indexa}{\PP V} :=
      \inv{\iota_{\indexa}}\Bigl(\sigma_r \left(\PP(S^{\indexa} V) \times \PP(S^{d- \indexa} V)\right)\Bigr).
\end{equation}
The defining equations of secant varieties of the Segre product of two projective spaces are well known
  --- these are just the $(r+1)\times(r+1)$ minors of the matrix, whose entries are linear coordinates on the tensor product.
The pullback by $\iota_{\indexa}$ of the matrix is called the \emph{$\indexa$-th catalecticant matrix}.
Thus $\Ranklocus{r}{\indexa}{d-\indexa}{\PP V}$
  is cut out by the $(r+1)\times(r+1)$ minors of the catalecticant matrix.

The image $\iota_{\indexa}(v_d(\PP V))$ is contained in $\PP(S^{\indexa} V) \times \PP(S^{d- \indexa} V)$
  and thus we always have
\begin{equation}\label{equ_embedding_of_sigma_for_Veronese}
  \sigma_r\bigl(v_d(\PP V)\bigr) \subset  \Ranklocus{r}{\indexa}{d-\indexa}{\PP V}.
\end{equation}
In general, $\Ranklocus{r}{\indexa}{d-\indexa}{\PP V}$ can be strictly larger than $\sigma_r\bigl(v_d(\PP V)\bigr)$,
even if we replace the right side of the inclusion \eqref{equ_embedding_of_sigma_for_Veronese} by the intersection of $\Ranklocus{r}{\indexa}{d-\indexa}{\PP V}$ for all $\indexa$.
The known examples are when $d$ is relatively small with respect to $r$ and $\dim V$.

The purpose of this article is to understand explicitly the locus $\Ranklocus{r}{\indexa}{d-\indexa}{\PP V}$
  for $d$ large and $\indexa$ sufficiently close to $\frac{d}{2}$.
In particular, Theorem~\ref{thm_secants_are_cut_out_by_catalecticants} proves
  $\sigma_r\bigl(v_d(\PP V)\bigr) =  \Ranklocus{r}{\indexa}{d-\indexa}{\PP V}$ as sets,
  whenever either  $\dim V$ or $r$ is small.
Eisenbud asked when does the equality hold~ \cite[Question~1.2.2]{jabu_ginensky_landsberg_Eisenbuds_conjecture}.
In general, the secant variety is strictly smaller:

\begin{thm}\label{thm_counterexamples_to_Eisenbud}
   Suppose either
   \begin{itemize}
     \item $n \ge 6$ and $r \ge 14$ or
     \item $n = 5$ and $r \ge 42$ or
     \item $n = 4$ and $r \ge 140$.
   \end{itemize}
   If $d \ge 2r-1$ and $V \simeq \CC^{n+1}$, then set-theoretically
   \[
      \sigma_r\bigl(v_d(\PP V)\bigr)  \varsubsetneq \bigcap_{\indexa \in \setfromto{1}{d-1}} \Ranklocus{r}{\indexa}{d-\indexa}{\PP V}.
   \]
\end{thm}

We prove this theorem in Section~\ref{section_smoothability_of_Gorenstein}.

\begin{rem}
  We learned from Anthony Iarrobino after the submission, that Theorem~\ref{thm_counterexamples_to_Eisenbud}
     has a significant overlap with \cite[Cor.~6.36]{iarrobino_kanev_book_Gorenstein_algebras}.
  In particular, the answer to the question of Eisenbud was previously known to be negative in many cases.
  In brief, \cite[Cor.~6.36]{iarrobino_kanev_book_Gorenstein_algebras} has better bounds on $d$,
     but (except in case $n=5$) worse bounds on $r$.
  Also some statements of Lemmas~\ref{lemma_hilbert_function_of_Omega_p_is_nice} and \ref{lemma_J_is_saturated} are in
     \cite[Prop.~C.33]{iarrobino_kanev_book_Gorenstein_algebras}.
\end{rem}

\subsection{Secant variety versus cactus variety}

We introduce \emph{the $r$-th cactus variety} of $X \subset \PP W$,
  which we denote $\cactus{r}{X}$; see Section~\ref{section_cactus} for details
  and Appendix~\ref{section_why_cactus} for explanation of the name.
The cactus variety $\cactus{r}{X}$ is the closure of union of
  the scheme-theoretic linear spans of $R$,
  where $R$ runs through all zero-dimensional subschemes of $X$ of degree at most $r$
  --- see also \eqref{equ_define_cactus2}.
This variety fits in between the secant and the zero locus of the minors
\[
  \sigma_r\bigl(v_d(\PP V)\bigr) \subset \cactus{r}{v_d(\PP V)} \subset \Ranklocus{r}{\indexa}{d-\indexa}{\PP V}
\]
(see Propositions~\ref{prop_basic_properties_of_cactus} and \ref{prop_cactus_is_contained_in_rank_locus}).
Thus we split the study of the inclusion \eqref{equ_embedding_of_sigma_for_Veronese}
  into two steps.
The following theorem explains the relation between
  the equality $\sigma_r\bigl(v_d(X)\bigr) = \cactus{r}{v_d(X)}$
  and smoothability of Gorenstein zero-dimensional schemes
  (see Section~\ref{section_cactus} for definition of smoothability
    and Sections~\ref{section_flattenings_and_Gorenstein_Artin_algebras}
    and \ref{section_smoothability_of_Gorenstein}
    for an overview of the Gorenstein schemes and their smoothability).
\begin{thm}\label{thm_when_secant_equals_cactus}
   Suppose $X \subset \PP V$ is a projective variety and let $r \ge 1$ be an integer.
   We say that \ref{item_condition_on_smoothability_of_all_Gorenstein} holds if
\renewcommand{\theenumi}{$(\star)$}
   \begin{enumerate}
     \item \label{item_condition_on_smoothability_of_all_Gorenstein}
         every zero-dimensional Gorenstein subscheme of $X$ of degree at most $r$ is smoothable in $X$.
   \end{enumerate}
   Then
\renewcommand{\theenumi}{(\roman{enumi})}
   \begin{enumerate}
      \item \label{item_if_star_holds_then_sigma_eq_cactus}
         If \ref{item_condition_on_smoothability_of_all_Gorenstein} holds, then $\sigma_r(X) = \cactus{r}{X}$.
      \item \label{item_if_sigma_eq_cactus_then_star_holds}
         If $\sigma_r\bigl(v_d(X)\bigr) = \cactus{r}{v_d(X)}$ for some $d \ge 2r-1$,
            then \ref{item_condition_on_smoothability_of_all_Gorenstein} holds.
   \end{enumerate}
\end{thm}
We prove this theorem in Section~\ref{section_cactus}.

Note that \ref{item_condition_on_smoothability_of_all_Gorenstein}
  is independent of the embedding of $X$.
Thus in the situation of \ref{item_if_star_holds_then_sigma_eq_cactus},
  also $\sigma_r\bigl(v_d(X)\bigr) = \cactus{r}{v_d(X)}$  for all $d$.

Let  $\HilbGor{r}{X}$ be  the subset of the Hilbert scheme of $r$ points
  parametrising Gorenstein schemes.
For an irreducible $X$ with $\dim X \ge 1$
  the condition \ref{item_condition_on_smoothability_of_all_Gorenstein}
  is equivalent to irreducibility of  $\HilbGor{r}{X}$.
The latter  is intensively studied,
  see Section~\ref{section_smoothability_of_Gorenstein} for an overview and references.
In particular,
   the condition \ref{item_condition_on_smoothability_of_all_Gorenstein}
   is known to hold for smooth $X$ if $\dim X \le 3$, or if $r\le 10$
   (see Proposition~\ref{prop_when_Gorenstein_are_smoothable}).
On the other hand, it is known to fail if $\dim X \ge 6$ and $r \ge 14$
   and it fails for $\dim X = 4,5$, for $r$ sufficiently large
   (see Proposition~\ref{prop_when_Gorenstein_are_not_smoothable}).
It also fails for many singular $X$.
For instance, consider a double point $R \simeq \Spec \CC[x]/ \langle x^2\rangle$.
Then embeddings $R \subset X$ with support at $x\in X$
   correspond to points in the (projectivised) Zariski tangent space at $x$,
   whereas $R \subset X$ is smoothable in $R$
   if and only if the tangent direction of $R$ is contained in the \emph{tangent star} of $X$ at $x$
   (see \cite[\S1.4]{jabu_ginensky_landsberg_Eisenbuds_conjecture}).
The dimension of the tangent star is at most $2\dim X$,
   so it is easy to construct examples where the tangent star is smaller than the Zariski tangent space ---
   for instance any curve with a singularity, which is not isomorphic to a planar singularity will do.
These properties are exploited in \cite[\S3]{jabu_ginensky_landsberg_Eisenbuds_conjecture}.

\subsection{Cactus variety versus catalecticant minors}

The next step is to understand when
  $\cactus{r}{v_d(\PP V)} = \Ranklocus{r}{\indexa}{d-\indexa}{\PP V}$ as sets.
We claim that for $d$ sufficiently large this equality always holds.

\begin{thm}\label{thm_our_version_of_Eisenbuds_for_Pn}
   Suppose $d\ge 2r$ and $r \le \indexa  \le d-r$.
   Then $\cactus{r}{v_d(\PP V)} = \Ranklocus{r}{\indexa}{d-\indexa}{\PP V}$ as sets.
\end{thm}

The proof of this theorem is effective in the sense,
  that given a point in $[p] \in \Ranklocus{r}{\indexa}{d-\indexa}{\PP V}$,
  we can explicitly find the unique
    (see Theorem~\ref{thm_our_version_of_Eisenbuds_for_Pn_effective_version})
  smallest scheme $R \subset \PP V$,
  such that $[p] \in \langle v_d(R) \rangle$.
For this purpose
  let $\annihp \subset \Sym V^*$  be the \emph{annihilator} of $p \in S^d V$,
  that is after the identification of $\Sym V^*$
  with the algebra of polynomial differential operators with constant coefficients,
  $\annihp$ is the ideal of the operators annihilating $p$
  (see Section~\ref{section_flattenings_and_Gorenstein_Artin_algebras} for further details).

\begin{thm}\label{thm_our_version_of_Eisenbuds_for_Pn_effective_version}
   Suppose $d\ge 2r$ and $r \le \indexa  \le d-r$.
   Let $p\in S^d V$ be such that $[p] \in \Ranklocus{r}{\indexa}{d-\indexa}{\PP V}$.
   Let $\ccJ$ be the homogeneous ideal generated by the first $\indexc$ degrees of $\annihp$,
     where $\indexc$ is any number such that $r \le \indexc \le d-r+1$.
   Let $R \subset \PP V$ be the scheme defined by $\ccJ$.
   Then
     \begin{enumerate}
        \item \label{item_R_has_the_right_dim_deg_and_has_p_in_span}
              $\dim R = 0$, $\deg R \le r$, and $[p] \in \langle v_d (R) \rangle$;
        \item \label{item_R_is_unique}
              $R$ is the smallest with respect to inclusion:
                 if $Q \subset \PP V$ is another scheme with
                 $\dim Q = 0$, $\deg Q \le d-r+1$, and $[p] \in \langle v_d (Q) \rangle$,
                 then $R \subset Q$;
              In particular, $R$ is unique such scheme of minimal degree.
        \item \label{item_J_is_saturated_and_independent_of_choice}
              $\ccJ$ is a saturated ideal, independent of the choice of $\indexc$.
     \end{enumerate}
\end{thm}

This has some interesting consequences, especially
  if $p$ is on an honest secant $\PP^{r-1}$,
  as we are able to determine the linear forms $v_i \in V$, giving
  $p = {v_1}^d + \dotsb + {v_r}^d$.

\begin{cor}\label{cor_decomposition_for_points_on_pure_secants}
   Suppose $d$, $r$, $\indexa$, $p$ and $R$ are as in Theorem~\ref{thm_our_version_of_Eisenbuds_for_Pn_effective_version}.
   Then $p = {v_1}^d + \dotsb + {v_r}^d$ for some $v_i \in V$
      if and only if $R$ is reduced.
   Moreover, in such a case, $R  \subset \setfromto{[v_1]}{[v_r]}$
      with equality if $p \notin \Ranklocus{r-1}{\indexa}{d-\indexa}{\PP V}$.
\end{cor}

We prove Theorems~\ref{thm_our_version_of_Eisenbuds_for_Pn}
              and \ref{thm_our_version_of_Eisenbuds_for_Pn_effective_version},
   and  Corollary~\ref{cor_decomposition_for_points_on_pure_secants}
   in Section~\ref{section_bounds_on_Hilb_functions}.

\subsection{Secant varieties of Veronese reembeddings}

We also generalise our results to the following setup, which is motivated by
  the question of Eisenbud ---
  see \cite[Question~1.2.3]{jabu_ginensky_landsberg_Eisenbuds_conjecture}
  or Question~\ref{conj_restricted_Eisenbud} below.
For a projective (possibly reducible) variety $X \subset \PP V$ let
\[
  \Ranklocus{r}{\indexa}{d-\indexa}{X} := \Ranklocus{r}{\indexa}{d-\indexa}{\PP V} \cap \langle v_d(X) \rangle,
\]
where $\langle v_d(X) \rangle$ is the linear span of $v_d(X) \subset \PP S^d V$.
Thus $\Ranklocus{r}{\indexa}{d-\indexa}{X}$ is defined by $(r+1)\times(r+1)$ minors of
  the catalecticant matrix with some linear substitutions determined by the embedding of $X$ into the projective space.

\begin{thm}\label{thm_our_version_of_Eisenbuds_for_X}
   Let $X \subset \PP V$ and $r \ge 1$.
   There exists an integer $d_0  = d_0(r,X)$
     such that for all $d \ge d_0$ and $r \le \indexa \le d-r$ the following equality of sets holds
   \[
     \cactus{r}{v_d(X)} = \Ranklocus{r}{\indexa}{d-\indexa}{X}.
   \]
   Here $d_0 = \max\set{2r, Got(h_X) + r - 1}$ and $Got(h_X)$
     is the Gotzmann number of the Hilbert polynomial of $X$
     (see \cite[Prop.~2.1.2]{jabu_ginensky_landsberg_Eisenbuds_conjecture}).
\end{thm}

We conclude that if, for instance, $X$ is smooth and $\dim X\le 3$ or $r \le 10$,
  then the answer to the set-theoretic version of Eisenbud's question is positive.

\begin{cor}\label{cor_Eisenbuds_holds_for_X}
   Let $X \subset \PP V$ and $r \ge 1$.
     First suppose condition \ref{item_condition_on_smoothability_of_all_Gorenstein}
     of Theorem~\ref{thm_when_secant_equals_cactus} holds.
   Then for all $d \ge d_0$ and $r \le \indexa \le d-r$ the following equality of sets holds:
   \[
     \sigma_{r}{(v_d(X))} = \Ranklocus{r}{\indexa}{d-\indexa}{X}.
   \]
   On the other hand, if \ref{item_condition_on_smoothability_of_all_Gorenstein} fails to hold,
     then for all $d \ge 2r-1$ set-theoretically
   \[
     \sigma_{r}{(v_d(X))} \varsubsetneq \bigcap_{\indexa=1}^{d-1} \Ranklocus{r}{\indexa}{d-\indexa}{X}.
   \]
\end{cor}

Both the theorem and its corollary are proved in Section~\ref{section_arbitrary_X}.

\subsection*{Overview}

In Section~\ref{section_cactus} we introduce the cactus variety
  and compare it with the secant variety.
In Section~\ref{section_flattenings_and_Gorenstein_Artin_algebras}
  we begin the comparison of cactus variety with the locus determined by catalecticant minors
  and prove the easier inclusion.
In Section~\ref{section_geometry_of_gorenstein} we motivate the line of our argument
  in Section~\ref{section_bounds_on_Hilb_functions}, where
  we use Macaulay's bound on growth of Hilbert function and Gotzmann's Persistence Theorem
  to prove the other inclusion.
This is the most technical part of the article.
In Section~\ref{section_smoothability_of_Gorenstein}
  we review what is known about smoothability of Gorenstein schemes
  and conclude with the proofs of Theorems~\ref{thm_secants_are_cut_out_by_catalecticants}
  and \ref{thm_counterexamples_to_Eisenbud}.
In Section~\ref{section_arbitrary_X} we briefly review the history related to the Eisenbud's question
  and provide the generalisations of our results for projective space to arbitrary projective variety.
In Section~\ref{section_improving_bounds} we explain how to slightly improve the bounds on integers in our theorems,
  but relying only on an unpublished work in progress or on proofs which we only sketch.
In Appendix~\ref{section_why_cactus} we explain the similarities of the cactus variety
 and the plant of the family Cactaceae.

\subsection*{Acknowledgements}

The authors would like to thank Joseph Landsberg for suggesting this topic for research
  and for all the related discussions,
  Frank Schreyer, Gianfranco Casnati, Roberto Notari, Vivek Shende, Grzegorz Kapustka, Micha\l{} Kapustka and Gavin Brown
  for their hints and suggestions about Gorenstein schemes
  and Laurent Manivel for listening and his comments.
The authors thank the organisers of 2010 IMPANGA school at B\c{e}dlewo,
  where the authors (among many other scientific and non-scientific attractions) could meet F.~Schreyer, G. and M.~Kapustka.
The second author thanks Gianfranco Casnati, for invitation to Politecnico di Torino.
We are sincerely grateful the anonymous referee for his numerous suggestions how to improve the presentation.
Finally, we dedicate this article to our son Mi\l{}osz, who
  was very kind to us and healthy enough to allow a successful work on the problem
  --- both before, and after he was born.

\section{Cactus variety}\label{section_cactus}

For a scheme $R \subset \PP W$ by $\langle R \rangle$ we denote its scheme theoretic linear span,
that is the smallest linear subspace $\PP^q \subset \PP W$ containing $R$.

For a projective variety $X \subset \PP V$,
 we denote by $\Hilb_r (X)$ the Hilbert scheme of subschemes of $X$ of dimension $0$ and degree $r$ and
by $\Hilb_{\le r} (X)$ we denote the Hilbert scheme of subschemes of $X$ of dimension $0$ and degree at most $r$, that is:
\[
  \Hilb_{\le r} (X) = \bigsqcup_{q=1}^{r} \Hilb_{q} (X).
\]
Thus closed points of $\Hilb_r (X)$ are in one-to-one correspondence with subschemes of $X$ of dimension $0$ and degree $r$
  and by a slight abuse of notation we will write $R \in \Hilb_r (X)$ to mean that $R \subset X$ is the corresponding subscheme
  and vice versa.

Let $X$ be a projective variety.
We say that a zero-dimensional subscheme $R\subset X$ of degree $r$ is \emph{smoothable in $X$},
  if it is a flat limit of $r$ distinct points on $X$.
We say $R$ is \emph{smoothable} if it is smoothable in some smooth projective variety $X$.
In fact, if $R$ is smoothable, then it is smoothable in any smooth $X$:
\begin{prop}\label{prop_smothable_in_X_iff_smoothable_in_Y}
  Suppose $R$ is a zero-dimensional scheme of finite length $r$  and $X$ and $Y$ are two projective varieties.
  If $R$ can be embedded in $X$ and in $Y$, and $R$ is smoothable in $Y$,
    and $R \subset X$ is supported in the smooth locus of $X$, then $R$ is smoothable in $X$.
\end{prop}
Despite the proposition is a standard folklore fact according to experts,
  the authors were not able to find an explicitly identical statement in the literature.
The most similar statement is \cite[Lem.~2.2]{casnati_notari_irreducibility_Gorenstein_degree_9},
  where $X \simeq \PP W$ and $Y \simeq \PP W'$.
%  (note that the content in \cite{casnati_notari_irreducibility_Gorenstein_degree_9}
%   is not identical to the content of the arxiv version of that article).
Other related statements are \cite[Lem.~4.1]{cartwright_erman_velasco_viray_Hilb8},
  \cite[p.4]{artin_deform_of_sings} or \cite[p.2]{erman_velasco_syzygetic_smoothability}.
Below we present a straightforward reduction of the general case
  to the case of \cite[Lem.~2.2]{casnati_notari_irreducibility_Gorenstein_degree_9}.

\begin{prf}
Let $Y \subset \PP W'$ be the embedding of the projective variety $Y$. If $R$ is smoothable in $Y$,
    then it is automatically also smoothable in its ambient projective space $\PP W'$.
 Moreover, by \cite[Lem.~2.2]{casnati_notari_irreducibility_Gorenstein_degree_9}
    it is smoothable in any other embedding into projective space.

  It is enough to prove every irreducible component is smoothable in $X$,
    thus for simplicity we assume $R$ is supported at a single point.
Assume this point is $x \in X$, when considering the embedding $R \subset X$.
  We pick a small open analytic neighbourhood $D^n \subset X$ of $x\in X$
    and a holomorphic embedding $\phi\colon D^n \hookrightarrow \PP W$,
    where $\dim \PP W = \dim X$.
  Denote $p:= \phi(x) \in \PP W$.
  Thus $\phi|_{R} \colon R \to \phi(R)$ is an isomorphism of abstract schemes,
     and by our assumptions $\phi(R)$ is smoothable in $\PP W$
     and $\phi(R)$ is supported at $p$.
  Therefore there exist a curve $C$ with a point $c \in C$
    and  $ \Gamma \subset \PP W \times C$
    such that the projection $\Gamma \to C$ is flat,
    the general fibre is a union of $r$ reduced points,
    and the special fibre $\Gamma_c \subset \PP W \times \set{c}$ is equal to $\phi(R)$.
  The preimage $\tilde{\Gamma}:=\inv{(\phi\times \id_C)}(\Gamma)$ is therefore a \emph{holomorphic} smoothing of $R$ in $X$.
  Let $C_0 \subset C$ be a small analytically open neighbourhood of $c$.
  As we argue below, there is an induced holomorphic map $\xi_{\tilde{\Gamma}}\colon C_0 \to \Hilb_r(X)$,
     which maps $x\mapsto R$ and a general point to $r$ distinct points.
  Therefore $R$ is in the same irreducible component of $\Hilb_r(X)$ as $r$ distinct points,
     which proves the claim of proposition. It remains to explain the existence of $\xi_{\tilde{\Gamma}}$.

  Since $\tilde{\Gamma}$ is not an algebraic family,
    the existence of $\xi_{\tilde{\Gamma}}$ is not guaranteed directly by the universal property of $\Hilb_r(X)$.
  However, it is enough to recall that $\Hilb_r(X)$ is constructed as a closed subset of a Grassmannian $Gr$,
    and the universal property of $\Hilb_r(X)$ is a restriction of the universal property of the Grassmannian
     (see for instance \cite[Thm~VI.22 and pp.263--264]{eisenbud_harris}).
  Furthermore, the universal property of Grassmannian is the same in both algebraic and analytic categories.
  Thus there exists a map $\xi_{\tilde{\Gamma}}\colon C_0 \to Gr$
    (roughly, mapping a point $b \in C_0$ to the scheme-theoretic linear span of the $\tilde{\Gamma}_b$
     under some high degree Veronese reembedding).
  Moreover, the image of $\xi_{\tilde{\Gamma}}$ is contained in $\Hilb_r(X) \subset Gr$.
  Thus the claim is proved.
%THINK: Eisenbuda ksiazka, Corollary 7.5
\end{prf}

The relations of smoothable schemes and secant varieties are exploited in particular by
  \cite{bernardi_gimigliano_ida} and \cite{jabu_ginensky_landsberg_Eisenbuds_conjecture}.

Suppose $X \subset \PP W$ is a projective variety.
The $r$-th secant variety of $X$ is defined as the closure of the union of linear spans of
  $r$ distinct points on $X$:
\begin{align}
   \sigma_{r}(X) &= \overline{\bigcup \set { \langle R \rangle \mid R\subset X,  \ R=\text{at most $r$ distinct points in X}}}.
   \nonumber\\
   \intertext{%
Since there is ``closure'' in this definition, we obtain the same object,
  if we add to the union the linear spans of the limiting schemes:}
  \sigma_{r}(X) &= \overline{ \bigcup \set { \langle R \rangle \mid   {R\in \Hilb_{\leq r}(X)}, R\text{ is smoothable in }X } }.
  \nonumber\\
  \intertext{%
Here we introduce a variant of the secant variety for a subvariety $X \subset \PP W$,
 and we call it the \emph{$r$-th cactus variety of $X$}:}
  \cactus{r}{X} &= \overline{ \bigcup \set { \langle R \rangle \mid {R\in \Hilb_{\leq r}(X)}} },\label{equ_define_cactus2}
\end{align}
where this time $R$ runs through \emph{all} zero-dimensional subschemes of degree at most $r$
(we forget the smoothability requirement).

Note however, that the expected dimension of $\cactus{r}{X}$ is in general much bigger
  than the expected secant dimension, because the dimension of some components
  of the Hilbert scheme $\Hilb_r X$ could be large.

\begin{prop}\label{prop_basic_properties_of_cactus}
   We have the following elementary properties:
   \begin{enumerate}
      \item \label{item_sigma_subset_cactus}
               $\sigma_r(X) \subset \cactus{r}{X} \subset \langle X \rangle$
               and $\cactus{q}{X} \subset \cactus{r}{X}$ for $q \le r$
               and $\cactus{r}{X} \subset \cactus{r}{Y}$ for $X \subset Y \subset \PP W$.
      \item \label{item_Gorenstein_are_enough_for_cactus}
               In the definition of cactus variety \eqref{equ_define_cactus2} it is enough to consider
               Gorenstein schemes:
               \[
                 \cactus{r}{X} = \overline{ \bigcup \set { \langle R \rangle \mid {R\in \Hilb_{\leq r}(X)},\text{ $R$ is Gorenstein}}}.
               \]
   \end{enumerate}
\end{prop}

\begin{prf}
   Part~\ref{item_sigma_subset_cactus} follows directly from the definitions.

   Part~\ref{item_Gorenstein_are_enough_for_cactus}
     follows from the following Lemma~\ref{lem_if_not_Gorenstein_then_covered}.
   The lemma implies, that the linear spans of non-Gorenstein schemes
in the definition of  $\cactus{r}{X}$ are redundant.
They are covered by linear spans of shorter schemes.
\end{prf}

\begin{lemma}\label{lem_if_not_Gorenstein_then_covered}
   If $R \subset \PP W$ is not Gorenstein, then
\[
        \langle R \rangle  = \bigcup \set{\langle Q \rangle \mid {Q \in \Hilb_{\le r-1} R}}.
\]
\end{lemma}
The lemma is very similar to
        \cite[Lemma 2.4.4(i) \& (iii)]{jabu_ginensky_landsberg_Eisenbuds_conjecture}.
   For the reader's convenience, we rewrite the detailed proof.
\begin{prf}
   Let $q = \deg R -1$.
   Since $R$ is not Gorenstein, $\dim \Hilb_{q} R  > 0$ by
     \cite[Lemma 2.4.4(i)\&(ii)]{jabu_ginensky_landsberg_Eisenbuds_conjecture}.
   Suppose $Q \subset R$ is a subscheme, such that $\deg Q = q$.
   Then either $ \langle Q \rangle = \langle R \rangle$ or $ \langle Q \rangle$
      is a hyperplane in $\langle R \rangle$.
   If the first case happens for some $Q$, then the conclusion of the lemma holds.
   Suppose only the second case happens.
   We claim $Q = \langle Q \rangle \cap R$, where the intersection is scheme-theoretic.
This is because $Q \subset \langle Q \rangle \cap R$, $\deg Q = q = \deg R -1$,
      so $ \deg R -1 \le \deg \bigl(\langle Q \rangle \cap R \bigr) \le \deg R$.
   But we cannot have $\deg \bigl(\langle Q \rangle \cap R \bigr) = \deg R$,
      as by our assumption $\langle Q \rangle \cap R \ne R$. So $\deg \bigl(\langle Q \rangle \cap R \bigr) = q$
      and $Q = \langle Q \rangle \cap R$.

   Thus the hyperplane $\langle Q \rangle \subset \langle R \rangle$
      determines $Q$ uniquely
      and two different schemes $Q, Q' \subset R$ with $\deg Q = \deg Q' = q$ have
      $ \langle Q \rangle \ne \langle Q' \rangle$.
   Therefore, the locus:
     \begin{equation}\label{equ_union_of_hyperplanes}
        \bigcup \set{\langle Q \rangle \mid {Q \in \Hilb_{q} R}}.
     \end{equation}
     is a union of pairwise different hyperplanes in $\langle R \rangle$
     parametrised by projective scheme $\Hilb_{q} R$ of positive dimension.
   Thus it is a closed subset of $\langle R \rangle$ of dimension at least $\dim \langle R \rangle$,
     so the union \eqref{equ_union_of_hyperplanes} is equal to $\langle R \rangle$ and the lemma is proved.
\end{prf}

\begin{rem}\label{rem_enough_to_take_exactly_deg_r_schemes}
  If $X$ has at least $r$ points (that is, either $\dim X \ge 1$ or $\dim X =0$ and $\deg X \ge r$,
    then we can replace $R\in \Hilb_{\le r}(X)$ with $R\in \Hilb_{r}(X)$
    in the definitions of $\sigma_{r}(X)$ and $\cactus{r}{X}$ (see \eqref{equ_define_cactus2}), as well as in
    Proposition~\ref{prop_basic_properties_of_cactus}\ref{item_Gorenstein_are_enough_for_cactus}
    and also in Proposition~\ref{properties_of_cactus_of_Veronese_reembeddings} below.
  This is because we can always add to $R$ a bunch of distinct points, if needed.
  Even in the case $\dim X =0$ and $\deg X =t < r$,
    we have $\sigma_{r}(X) = \sigma_{t}(X)$ and $\cactus{r}{X} = \cactus{t}{X}$,
    it is enough to consider $R\in \Hilb_{t}(X)$.
  However, for some of the proofs it is more convenient to use all schemes in $\Hilb_{\le r}(X)$.
\end{rem}

In the cases of interest in this paper we can also get rid of the closure in \eqref{equ_define_cactus2}:
\begin{prop}\label{properties_of_cactus_of_Veronese_reembeddings}
   Suppose $d \ge r-1$.
   Then
   \[
      \cactus{r}{v_d(X)} = { \bigcup \set { \langle v_d(R) \rangle \mid {R\in \Hilb_{\leq r}(X)},\text{ $R$ is Gorenstein}}}.
   \]
   Similarly,
   \[
      \sigma_r(v_d(X)) = { \bigcup \set { \langle v_d(R) \rangle \mid {R\in \Hilb_{\leq r}(X)},\text{ $R$ is smoothable}}}.
   \]
\end{prop}

  For the secant variety, the statement is
    \cite[Lemma~2.1.5]{jabu_ginensky_landsberg_Eisenbuds_conjecture}
    and it is essentially a corollary of \cite[Prop.~11]{bernardi_gimigliano_ida}
    --- note however \cite[Prop.~11]{bernardi_gimigliano_ida} has an unnecessary smoothness assumption.
  For the cactus variety, the argument is identical, taking in account Lemma~\ref{lem_if_not_Gorenstein_then_covered}.
  For the sake of completeness, we reproduce both arguments bellow.

\begin{prf}
  For simplicity, we assume that $X$ has at least $r$ points (see Remark~\ref{rem_enough_to_take_exactly_deg_r_schemes}).
  To prove the statement for the cactus variety, let $\operatorname{H}_r:= \Hilb_r(X)_{red}$,
    the reduced subscheme of the Hilbert scheme.
  To prove the statement for the secant variety,
    let $\operatorname{H}_r$ be the reduced subscheme of irreducible component containing the smoothable schemes
    (or, if $X$ is not irreducible, then the reduced union of all such components).

  Let $Gr:=Gr(\PP^{r-1}, \PP(S^d V))$ be the Grassmannian of linearly embedded $\PP^{r-1}$ in $\PP(S^d V)$.
  Let $\ccU$ be the universal bundle, that is $\ccU \subset Gr \times \PP(S^d V)$ and
    the fibre of $\ccU$ over $E \in Gr$ is the linear space $E \subset \PP(S^d V)$.

  Since $d \ge r-1$, for every $R\in \operatorname{H}_r$ we have $ \dim \langle v_d (R) \rangle = r-1$
    (see \cite[Lem.~2.1.3]{jabu_ginensky_landsberg_Eisenbuds_conjecture}).
  Thus we obtain a well defined regular map $\operatorname{H}_r \to Gr$,
    $R\mapsto \langle v_d (R) \rangle$.
  In fact, this is precisely the restriction of embedding used to construct the Hilbert scheme $\Hilb_r(X)$,
    see for instance \cite[Thm~VI.22 and pp.263--264]{eisenbud_harris}.
  Let $\ccU_{\operatorname{H}_{r}}$ be the pullback of $\ccU$ under this map.
  That is $\ccU_{\operatorname{H}_r}$ is a $\PP^{r-1}$-bundle over $\operatorname{H}_r$ with a natural regular map
    $\pi \colon \ccU_{\operatorname{H}_r} \to \PP(S^d V)$,
    such that the fibre over $R \in \operatorname{H}_r$ is mapped onto $\langle v_d (R) \rangle \subset \PP(S^d V)$.
  By definition, the cactus variety $\cactus{r}{v_d(X)}$ or the secant variety $\sigma_r(v_d(X))$, respectively, is equal to $\overline{\pi(\ccU_{\operatorname{H}_r})}$.
  Both $\operatorname{H}_r$ and $\ccU_{\operatorname{H}_r}$ are projective, so the image of $\pi $ is closed,
    and thus $\pi(\ccU_{\operatorname{H}_r})$ is equal to either $\cactus{r}{v_d(X)}$
    or $\sigma_r(v_d(X))$, respectively.
  That is, for every $ p \in \cactus{r}{v_d(X)}$ or $ p \in\sigma_r(v_d(X))$, respectively,
    there exists $R \in \operatorname{H}_r$,
    such that $p \in \langle v_d (R) \rangle$, as claimed.
  Moreover, if $p\in \cactus{r}{v_d(X)}$, then by Lemma~\ref{lem_if_not_Gorenstein_then_covered},
    we may assume that $R$ is Gorenstein.
\end{prf}

We conclude with the proof of Theorem~\ref{thm_when_secant_equals_cactus}
  which claims that the equality between cactus variety and secant variety is equivalent
  to smoothability of Gorenstein schemes.

\begin{prf}[ of Theorem~\ref{thm_when_secant_equals_cactus}]
   Let $X \subset \PP V$ be a projective variety and $r\ge 1$ an integer.
   First we prove part~\ref{item_if_star_holds_then_sigma_eq_cactus}, so suppose
      all zero-dimensional Gorenstein schemes of degree at most $r$ are smoothable in $X$.
   By Proposition~\ref{prop_basic_properties_of_cactus}\ref{item_sigma_subset_cactus}
      the inclusion  $\sigma_r(X) \subset \cactus{r}{X}$ holds.
   To prove the other inclusion, let
   \[
      U := { \bigcup \set { \langle R \rangle \mid {R\in \Hilb_{\leq r}(X)},\text{ $R$ is Gorenstein}}}.
   \]
   By Proposition~\ref{prop_basic_properties_of_cactus}\ref{item_Gorenstein_are_enough_for_cactus}
       the set $U$ is dense in $\cactus{r}{X}$.
   But our assumption on smoothability of Gorenstein schemes implies that $U \subset \sigma_r (X)$.
   Since $\sigma_r (X)$ is closed, it follows that $\cactus{r}{X} \subset \sigma_r (X)$.

   To prove part~\ref{item_if_sigma_eq_cactus_then_star_holds},
       suppose $d \ge 2r-1$ and $\sigma_r(v_d(X)) = \cactus{r}{v_d(X)}$.
   Let $R\subset X$ be a Gorenstein zero-dimensional subscheme of degree at most $r$.
   We have to show that $R$ is smoothable in $X$.
   To prove the claim, first suppose $\deg R =r$.
   Take a general $p \in \langle v_d (R) \rangle \subset \cactus{r}{v_d(X)}$.
   Then by \cite[Lemma 2.4.4(iii)]{jabu_ginensky_landsberg_Eisenbuds_conjecture},
      $p \notin \langle v_d (R') \rangle$ for any non-trivial subscheme $R' \subsetneq R$,
      so that $R$ is \emph{minimal} (with respect to inclusion) such that $p \in \langle v_d (R)\rangle$.
   With our assumptions about $r$, $d$ and minimality of $R$, such $R$ is also \emph{unique}, in the sense
      it is the only zero-dimensional subscheme of $X$ of degree at most $r$,
      with $p \in \langle v_d (R) \rangle$ (see \cite[Cor.~2.2.1]{jabu_ginensky_landsberg_Eisenbuds_conjecture}).
   By our assumption $p \in \sigma_r(v_d(X))$,
      thus by Proposition~\ref{properties_of_cactus_of_Veronese_reembeddings}
      there exists $Q \subset X$, a zero-dimensional subscheme of degree at most $r$,
      which is smoothable in $X$ and $p \in \langle v_d (Q) \rangle$.
   Thus by uniqueness of $R$,
     we must have $R=Q$, so $R$ is smoothable in $X$ as claimed.

   If $\deg R <r$, then we can replace $R$ in the above considerations with
     $R':=R \cup \setfromto{x_{\deg R+1}}{x_r}$,
     where the $x_i\in X$ are some distinct reduced points, disjoint from $R$.
   This is always possible, unless $\dim X =0$ and $X$ has too few points,
     but then the claim trivially holds.
\end{prf}

\section{Flattenings and Gorenstein Artin algebras}\label{section_flattenings_and_Gorenstein_Artin_algebras}

For an introduction to the topic of Gorenstein Artin algebras see \cite[\S21.2]{eisenbud} and
\cite{iarrobino_kanev_book_Gorenstein_algebras}.
In Section~\ref{section_intro} we denoted by
$
  \iota_{\indexa}\colon \PP S^d V  \to \PP(S^{\indexa} V \otimes S^{d- \indexa} V)
$
the natural linear embedding.
By a slight abuse of notation we will use the same letter
  $\iota_{\indexa}$ to denote the underlying map of vector spaces:
\[
  \iota_{\indexa}\colon  S^d V  \to S^{\indexa} V \otimes S^{d- \indexa} V \simeq \Hom(S^{\indexa} V^*, S^{d- \indexa} V).
\]
For consistence, if $\indexa < 0$ or $\indexa > d$, then we assume $\iota_{\indexa}$ is identically $0$.

For $p \in S^d V$ and $\alpha \in S^{\indexa} V^*$ we denote by
  $\alpha \hook p:= \iota_{\indexa}(p)(\alpha)$.
This bilinear map $(\cdot \hook\cdot)\colon S^{\indexa} V^* \times S^{d} V \to S^{d-\indexa} V$
    (which, if interpreted as a tensor in $S^{\indexa} V \otimes S^{d} V^* \otimes S^{d-\indexa} V$,
      is the same tensor as $\iota_{\indexa}$)
  has various names in the literature:
  it is called \emph{flattening}, or \emph{contraction} and can be seen as a \emph{derivation}.
  That is, if we identify $\Sym V^*$ with the algebra of polynomial differential operators with constant coefficients,
    then $\alpha \hook p = \alpha (p)$, where on the right side $\alpha(p)$
    is seen as the application of the differential operator to $p$.
The following are natural properties of $\hook$, which we exploit in this section.
\begin{itemize}
 \item  For $\alpha \in S^{\indexa} V^*$, $\beta \in S^{\indexb} V^*$, $p \in S^d V$ we have
         \[
           \alpha \hook(\beta \hook p) = (\alpha \beta) \hook p = \beta \hook(\alpha \hook p).
         \]
 \item  If $\indexa = d$,  $\alpha \in S^{d} V^*$,  $p \in S^d V$,
          then $\alpha \hook p \in S^0 V \simeq \CC$ agrees with the natural pairing
          $(\alpha, p)$ arising from the duality
          $(\cdot, \cdot)\colon (S^{d} V)^* \otimes S^{d} V \to \CC$ (up to scale).
\end{itemize}

\begin{rem}
   The second item above seems to be slightly controversial.
   We received suggestions, both that it is trivial and not worth mentioning,
     as well as that it is false as stated.
   We illustrate the problem on the case $\dim V =2$ and $d=2$.
   Say a basis of $V$ is $(x,y)$ and the dual basis of $V^*$ is $(\alpha, \beta)$.
   Then $S^2 V$ has a basis $(x^2, xy, y^2)$,
      and $\alpha^2 \hook x^2 = 2$, $\alpha\beta \hook xy  = 1$ and $\beta^2 \hook y^2 =2$
      with the other products equal to zero.
   However, one could wrongly expect
      the natural pairing arising from the duality to be such that the monomial bases
      $(x^2, xy, y^2)$ and $(\alpha^2, \alpha \beta, \beta^2)$ are dual to each other
      (so that the products are all $0$ or $1$).
   This is however not true.
   If they were dual bases, then the pairing would depend on the choice of basis of $V$ (see below).
   The coefficients $2,1,2$
      (in general, these coefficients are the products of factorials of exponents of monomials)
      show up when we write the isomorphism $S^d (V^*) \simeq (S^d V)^*$ in coordinates.

   To argue that the second item is true, it is enough to observe,
      that both pairings are independent of the choices of coordinates, that is they are $GL(V)$-invariant.
   But the space of the $GL(V)$-invariant tensors in $S^d V^* \otimes S^d V$ is $1$-dimensional,
       which follows from the Pieri formula, see~\cite[Theorem 10.2.1]{procesi_book}.
\end{rem}

Suppose $p \in S^d V$.
The central objects of this section are the \emph{kernel} and the \emph{image} of $\cdot \hook p$.
So let $\annihp \subset \Sym V^*$ be the \emph{annihilator} of $p$,
  that is the homogeneous ideal defined by its homogeneous pieces:
\[
  \annihp^{\indexa}  := \set{ \alpha \in S^{\indexa} V^* \mid  \alpha \hook p = 0  } = \ker\bigl(\iota_{\indexa}(p)(\cdot)\bigr).
\]
with $\annihp := \bigoplus_{\indexa=0}^\infty \annihp^{\indexa} \subset \Sym V^*$.
We also define:
\[
  \Omega_p^{\indexa} := S^{\indexa} V^*  \hook p  = \im \bigl(\iota_{\indexa}(p)(\cdot)\bigr) \subset S^{d-\indexa} V.
\]
  and $\Omega_p:= \bigoplus_{\indexa=0}^d \Omega_p^{\indexa} \subset \Sym V$.

Note that naturally  $\Omega_p \simeq \Sym V^* / \annihp$ is a graded algebra,
  but the grading of $\Omega_p$ descending from $\Sym V^* / \annihp$ is reverse to
  the grading obtained from $ \Omega_p \subset \Sym V$.
Throughout the article we denote by $\Omega_p^{\indexa}$
  the $\indexa$-th graded piece with respect to the algebra grading.

\begin{rem}
   More generally, in the same way one can define $\annihp$ for non-homogeneous polynomial $p$,
     to obtain non-homogeneous ideal,   see \cite[Thm~21.6 and Exercise~21.7]{eisenbud} or
     \cite[\S2]{casnati_elias_notari_rossi_Poincare_series_and_deformations_of_Gorenstein_algebras_with_low_degree}.
   Although we are going to use this more general notion
     in the proof of Proposition~\ref{prop_when_Gorenstein_are_not_smoothable},
     here we only restrict our attention to the homogeneous case.
\end{rem}

\begin{rem}
   $\Omega_p$ seen as $\Sym V^* / \annihp$ is a Gorenstein Artin algebra.
   This is a classical construction, see \cite[Lem.~2.12]{iarrobino_kanev_book_Gorenstein_algebras}
      or \cite[Thm~21.6 and Exercise~21.7]{eisenbud} or
      \cite[\S2]{casnati_elias_notari_rossi_Poincare_series_and_deformations_of_Gorenstein_algebras_with_low_degree}
      and in fact, all such algebras arise in this way.
\end{rem}

Let $W \subset S^\indexa V^*$ be a linear subspace.
The map dual to this embedding is $S^{\indexa} V \twoheadrightarrow W^*$.
We denote the kernel of this map by $W^{\perp} \subset S^{\indexa} V$.

\begin{prop}\label{prop_basic_properties_of_Omega_p_and_I_p}
   Fundamental properties of $\annihp$ and $\Omega_p$ are:
   \begin{enumerate}
      \item \label{item_I_p_is_everything_for_high_delta}
            $\annihp^{\indexa} =  S^{\indexa} V^*$ for $\indexa >  d$, thus $\annihp$ defines an empty scheme in $\PP V$.
      \item \label{item_I_p_is_delta_saturated}
            For any $\indexb \le \indexa \le d$
	    we have
	    \[
	      \annihp^{\indexb}
	       = \set{\alpha \in S^{\indexb} V^* \mid S^{\indexa-\indexb} V^* \cdot \alpha  \subset \annihp^{\indexa}}.
	    \]
      \item \label{item_I_p_is_the_largest_ideal}
            $\annihp$ is the largest homogeneous ideal in $\Sym V^*$ such that $p \in (\annihp^{d})^{\perp}$,
            i.e.~if $\ccI$ is a homogeneous ideal such that  $p \in (\ccI^{d})^{\perp}$, then $\ccI \subset \annihp$.
      \item \label{item_Omega_p_contains_Omega_alpha_hook_p}
            For all $\alpha \in V^*$ we have $\Omega_{\alpha \hook p}^{\indexa-1} \subset \Omega_p^\indexa$.
      \item \label{item_Omega_is_perpendicular_to_I}
            $\Omega_p^{d-\indexa} =(\annihp^{\indexa})^{\perp} $,
              and $\Omega_p^{d-\indexa}$ is naturally dual to $\Omega_p^{\indexa}$.
      \item \label{item_S_delta_V_subset_I_perp}
            Suppose $\ccI \subset \Sym V^*$ is a homogeneous ideal such that $p \in (\ccI^d)^{\perp}$.
            Then for any $\indexa$, we have $\Omega_p^{d-\indexa} \subset (\ccI^{\indexa})^{\perp}$.
   \end{enumerate}
\end{prop}

\begin{prf}
   Part~\ref{item_I_p_is_everything_for_high_delta} is clear from the definitions.

   In \ref{item_I_p_is_delta_saturated} the inclusion $\subset$ is immediate, since $\annihp$ is an ideal.
   To prove $\supset$, let $\alpha \in S^{\indexb}V^*$ be such that
     $S^{\indexa- \indexb}V^* \cdot \alpha \subset \annihp^{\indexa}$.
   Then
   \[
      0=(S^{\indexa-\indexb}V^* \cdot \alpha) \hook p =   S^{\indexa-\indexb}V^* \hook (\alpha  \hook p).
   \]
   Thus $\alpha  \hook p$ is a homogeneous polynomial of degree $d-\indexb$,
     whose all degree $\indexa-\indexb$ derivatives are equal $0$.
   Since $(d-\indexb) - (\indexa-\indexb) \ge 0$,
     we must have $\alpha  \hook p =0$ and $\alpha \in \annihp^{\indexb}$.

   To prove \ref{item_I_p_is_the_largest_ideal},
     note that $\CC \cdot p = (\annihp^{d})^{\perp}$.
   Suppose $\ccI \subset \Sym V^*$ is a homogeneous ideal, such that $p \in (\ccI^d)^{\perp}$.
   Then $\ccI^d \subset \annihp^{d}$.
   Suppose $\alpha \in \ccI$ is a homogeneous polynomial.
   If $\deg \alpha > d$, then $\alpha \in \annihp$ by \ref{item_I_p_is_everything_for_high_delta}.
   Suppose $\indexa:=\deg \alpha \le d$.
   Then $S^{d-\indexa}V^* \cdot \alpha \subset I^d \subset \annihp^{d}$.
   By \ref{item_I_p_is_delta_saturated} this implies $\alpha \in \annihp^{\indexa}$.
   Therefore $\ccI \subset \annihp$.

  Part \ref{item_Omega_p_contains_Omega_alpha_hook_p} is clear from the definitions.

  For part \ref{item_Omega_is_perpendicular_to_I} consider
  the following dual short exact sequences:
  \begin{alignat*}{4}
     0 &\to \annihp^{\indexa} &&\to S^{\indexa} V^* &&\to \Omega_p^{\indexa} &&\to 0\\
     0 &\to (\Omega_p^{\indexa})^* &&\to S^{\indexa} V &&\to (\annihp^{\indexa})^* &&\to 0
  \end{alignat*}
  Here $ (\Omega_p^{\indexa})^* = (\annihp^{\indexa})^{\perp}$.
  Take any $f \in \Omega_p^{d-\indexa}$.
  We claim $f \in (\annihp^{\indexa})^{\perp}$, that is $\annihp^{\indexa} \hook f =0$.
  Since $f \in \Omega_p^{d-\indexa}$, there exists $\alpha \in S^{d-\indexa} V^*$
    such that $f = \alpha \hook p$. Hence:
  \[
     \annihp^{\indexa} \hook f = \annihp^{\indexa} \hook (\alpha \hook p) =  \alpha \hook(\annihp^{\indexa}  \hook p) =0
  \]
  Thus $\Omega_p^{d-\indexa} \subset (\Omega_p^{\indexa})^*$
    and in particular $\dim \Omega_p^{d-\indexa} \le \dim \Omega_p^{\indexa}$.
  If we exchange the roles of $\indexa$ and $d-\indexa$,
    we also get $\dim \Omega_p^{d-\indexa} \ge \dim \Omega_p^{\indexa}$,
  thus $\Omega_p^{\indexa} = (\Omega_p^{d-\indexa})^* = (\annihp^{\indexa})^{\perp}$ as claimed.

  To prove part~\ref{item_S_delta_V_subset_I_perp},
    we have $\ccI \subset \annihp$ by \ref{item_I_p_is_the_largest_ideal}.
  Thus by \ref{item_Omega_is_perpendicular_to_I}:
  \[
    \Omega_p^{d-\indexa} = (\annihp^{\indexa})^{\perp} \subset (\ccI^{\indexa})^{\perp}
  \]
\end{prf}

\begin{prop}\label{prop_if_p_in_span_of_R_then_Omega_in_span_of_R}
  Let $p \in S^d V$ and $R \subset \PP V$ be a scheme.
  If $[p] \in \langle v_d(R)\rangle$,
    then for all integers $i$ we have $\PP(\Omega_p^{d-\indexa}) \subset \langle v_{\indexa} (R) \rangle$.
\end{prop}

\begin{prf}
   Let $\ccJ \subset \Sym V^*$ be the saturated homogeneous ideal defining $R$.
   The condition $[p] \in \langle v_d (R) \rangle$ is equivalent to
    $p \in (\ccJ^d)^{\perp}$.
   By Proposition~\ref{prop_basic_properties_of_Omega_p_and_I_p}\ref{item_S_delta_V_subset_I_perp}:
   \[
     \PP (\Omega_p^{d-\indexa}) \subset \PP \left((\ccJ^{\indexa})^{\perp}\right) = \langle v_{\indexa} (R) \rangle.
   \]
\end{prf}

Recall, that $\Ranklocus{r}{\indexa}{d-\indexa}{\PP V} \subset \PP (S^d V)$
  is defined by vanishing of the $(r+1) \times (r+1)$ minors.
That is, for $p \in S^d V$,
  we have  $[p] \in \Ranklocus{r}{\indexa}{d-\indexa}{\PP V}$
  if and only if
  $\rk \bigl( (\cdot \hook p)\colon  S^{\indexa} V^* \to S^{d-\indexa} V  \bigr) \le r$.
Or, equivalently, $\dim \Omega_p^{\indexa} \le r$.

\begin{prop}\label{prop_cactus_is_contained_in_rank_locus}
   For all $r$, $d$ and $\indexa$, we have
   $\cactus{r}{v_d (\PP V)} \subset \Ranklocus{r}{\indexa}{d-\indexa}{\PP V}$.
\end{prop}

\begin{prf}
   Suppose $[p] \in \langle v_d (R) \rangle$ for a zero-dimensional subscheme $R \subset \PP V$ of degree at most $r$.
   We claim $[p] \in \Ranklocus{r}{\indexa}{d-\indexa}{\PP V}$.
   By Proposition~\ref{prop_if_p_in_span_of_R_then_Omega_in_span_of_R} we have
     $\PP (\Omega_p^{d-\indexa}) \subset \langle v_{\indexa} (R) \rangle$.
   Since the degree of $R$ is bounded by $r$,
     the dimension of $\langle v_{\indexa } (R) \rangle$ is bounded by $r-1$,
     and so  $\dim \Omega_p^{\indexa} \le r$.
   Thus $p \in \Ranklocus{r}{\indexa}{d-\indexa}{\PP V}$ and we obtain:
   \[
      \bigcup \set{\langle v_d( R )\rangle \mid {R\in \Hilb_{\leq r}(X)}} \subset \Ranklocus{r}{\indexa}{d-\indexa}{\PP V}.
   \]
   The closure of the left side (with reduced structure) is $\cactus{r}{v_d(\PP V)}$.
   The right side is closed,
     and thus we conclude $\cactus{r}{v_d(\PP V)} \subset \Ranklocus{r}{\indexa}{d-\indexa}{\PP V}$.
\end{prf}

\section{Geometry of Gorenstein Artin algebra}\label{section_geometry_of_gorenstein}

In the next section we prove that under suitable assumptions
  $\Ranklocus{r}{\indexa}{d-\indexa}{\PP V} \subset \cactus{r}{v_d(\PP V)}$,
  see Theorem~\ref{thm_our_version_of_Eisenbuds_for_Pn}.
In this section, we motivate our argument in Section~\ref{section_bounds_on_Hilb_functions}
  by presenting some geometric interpretations of
  statements in Proposition~\ref{prop_basic_properties_of_Omega_p_and_I_p}.
Formally, this section is not necessary for the overall argument.
However, the reader might benefit from reading this section
  by better understanding the later arguments.

We start with $[p] \in \Ranklocus{r}{\indexa}{d-\indexa}{\PP V}$,
  say $\dim \Omega_p^{\indexa} =r$,
  and we want to find a zero-dimensional subscheme $R \subset \PP V$ of degree $r$,
  such that $[p] \in \langle v_d (R) \rangle$.
In this section we address the following question:
\begin{itemize}
   \item[] What is the scheme $R$ that contains in its linear span given polynomial $[p] \in \langle v_d (R) \rangle$, where $p$ satisfies the above assumptions.

%With the above assumptions, suppose  $[p] \in \langle v_d (R) \rangle$; what is $R$?
\end{itemize}
Proposition~\ref{prop_if_p_in_span_of_R_then_Omega_in_span_of_R} gives a condition on $R$, without any additional assumptions. If, moreover, we have assumptions as in Theorem~\ref{thm_our_version_of_Eisenbuds_for_Pn},
then we can identify the unique candidate for the scheme $R$.

\begin{prop}
  With $p$, $d$, $R$ as in Proposition~\ref{prop_if_p_in_span_of_R_then_Omega_in_span_of_R},
  suppose in addition $\indexa \ge r$,  $\dim \Omega_p^{d-\indexa} =r$, $\dim R =0$ and $\deg R =r$.
  Then $v_{\indexa} (R) = \PP(\Omega_p^{d-\indexa}) \cap v_{\indexa} (\PP V)$.
\end{prop}

\begin{prf}
  Let $Q \subset \PP V$ be the scheme such that
     $v_{\indexa}(Q) = \PP(\Omega_p^{d-\indexa}) \cap v_{\indexa} (\PP V)$.
  The assumptions imply $\dim \langle v_{\indexa} (R) \rangle \le r-1 = \dim \PP (\Omega_p^{d-\indexa})$.
  By Proposition~\ref{prop_if_p_in_span_of_R_then_Omega_in_span_of_R},
     we have $\langle v_{\indexa} (R) \rangle  =  \PP (\Omega_p^{d-\indexa})$.
  Therefore $R \subset Q$.
  It only remains to prove $\dim Q =0$ and $\deg Q = r$.

  Suppose, by contradiction, there exists a subscheme $Q' \subset Q$, with $\dim Q' =0$ and $\deg Q' =r+1$.
  Since $\indexa \ge r$, by \cite[Lem.~2.1.3]{jabu_ginensky_landsberg_Eisenbuds_conjecture},
    the dimension of $ \langle v_{\indexa} (Q') \rangle$ is at least~$r$.
  This contradicts
    $\langle v_{\indexa} (Q') \rangle \subset \langle v_{\indexa} (Q) \rangle \subset \PP(\Omega_p^{\indexa})$
    and $\dim \Omega_p^{\indexa} =r$.
\end{prf}

Thus, if $R$ exists, as we claim in Theorems~\ref{thm_our_version_of_Eisenbuds_for_Pn}
  and \ref{thm_our_version_of_Eisenbuds_for_Pn_effective_version},
  then the ideal $\ccJ$ defining $R$ is generated by $\bigoplus_{\indexb=0}^{\indexa}\annihp^{\indexb}$.
The main problem is to show that $\ccJ$ defined in such a way is saturated.
In principle, $\ccJ$ could define an empty set in $\PP V$,
  and in fact this happens, if $\indexa$ is taken too large with respect to $r$ and $d$.
 More precisely, for certain $p$, for all $\indexa > d-r+1$, the ideal $\ccJ$ defines an empty set.
We pursue the problem  of~$\ccJ$ being saturated in the next section,
  by a careful examination of the Hilbert functions of~$\Omega_p$ and $\Psi := Sym V^* / \ccJ$.
The tools we use are Macaulay's bound on growth of Hilbert function of an algebra
  and Gotzmann's Persistence Theorem.

\section{Bounds on Hilbert functions of graded algebras}\label{section_bounds_on_Hilb_functions}

Macaulay \cite{macaulay_enumeration} obtained a bound on growth of the Hilbert function of graded algebras.
We use this bound in the form presented in \cite[Thm~3.3]{green_generic_initial_ideals} or
\cite[Thm~2.2(i), (iii)]{stanley_hilbert_functions_of_graded_algebras}.
In our situation, we use the symbol $h^{\langle i \rangle}$
  (defined in \cite[Def.~3.2]{green_generic_initial_ideals}
    or Equation (4) in \cite{stanley_hilbert_functions_of_graded_algebras})
  only when $ h \le i$,
  and in this case $h^{\langle i \rangle} =h$.
Thus the following statement is an immediate consequence of the Macaulay's bound:
\begin{cor}\label{cor_Macaulays_bound_for_the_growth_of_Hilbert_function}
   Let $\Psi$ be a graded algebra generated in degree $1$.
   Suppose for some $\indexa$, we have $\dim \Psi^{\indexa} \le \indexa$.
   Then  $\dim \Psi^{\indexa+1} \le \Psi^{\indexa}$.
\end{cor}

The following lemma is an easy consequence of this bound applied to $\Psi = \Omega_p$ for $p \in S^d V$.
It explains the Hilbert function of $\Omega_p$ behaves nicely
  (in particular, it is \emph{unimodal}, that is non-decreasing for values $<\frac{d}{2}$).

\begin{figure}[htb]
\centering
\textit{Figures~\ref{figure_with_dim_Omega_p_1}--\ref{figure_with_dim_Omega_p_4}.
        Illustrative proof of Lemma~\ref{lemma_hilbert_function_of_Omega_p_is_nice}}

\bigskip

\begin{tabular}{cc}
\begin{minipage}[t]{0.47\textwidth}
\includegraphics[width=\textwidth]{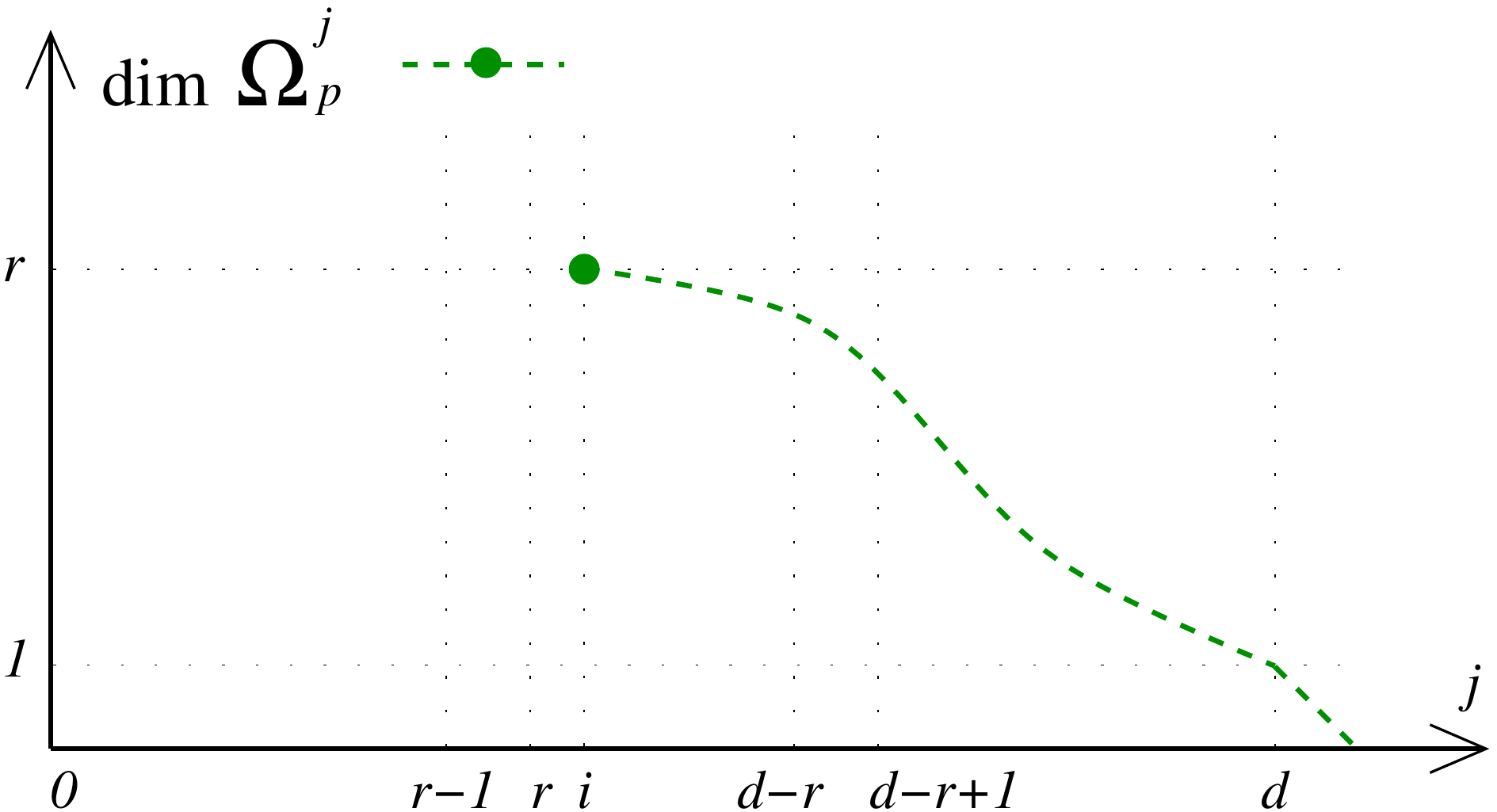}
\captionof{figure}{
We start with $\dim \Omega_p^{\indexa} = r$
% (marked by a $\bullet$ on the graph)
  and conclude from Corollary \ref{cor_Macaulays_bound_for_the_growth_of_Hilbert_function}
  that $\dim \Omega_p^{\indexb}$ is non-increasing for $\indexb \ge \indexa$.
Thus the last part of the lemma is proved.
}\label{figure_with_dim_Omega_p_1}
\end{minipage}
&
\begin{minipage}[t]{0.47\textwidth}
\includegraphics[width=\textwidth]{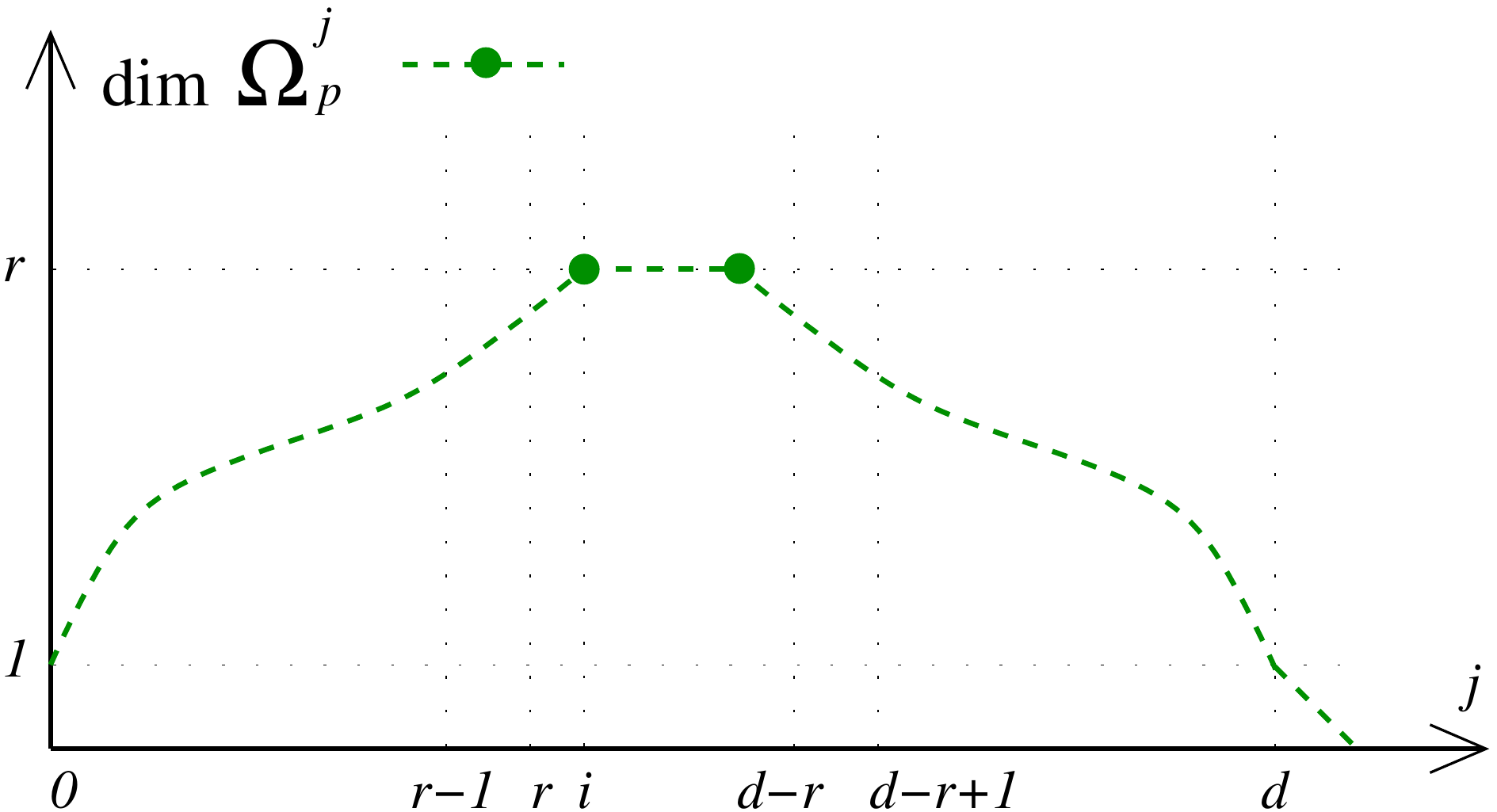}
\captionof{figure}{
We use the symmetry again to conclude $\dim \Omega_p^{\indexb}$ is non-decreasing for $\indexb \le d-\indexa$.
Thus the first part of the lemma is proved.
}\label{figure_with_dim_Omega_p_2}
\end{minipage}\\
\\
% \end{tabular}
% \end{figure}
% \begin{figure}[htb]
% \centering
% \begin{tabular}{cc}
\begin{minipage}[t]{0.47\textwidth}
\includegraphics[width=\textwidth]{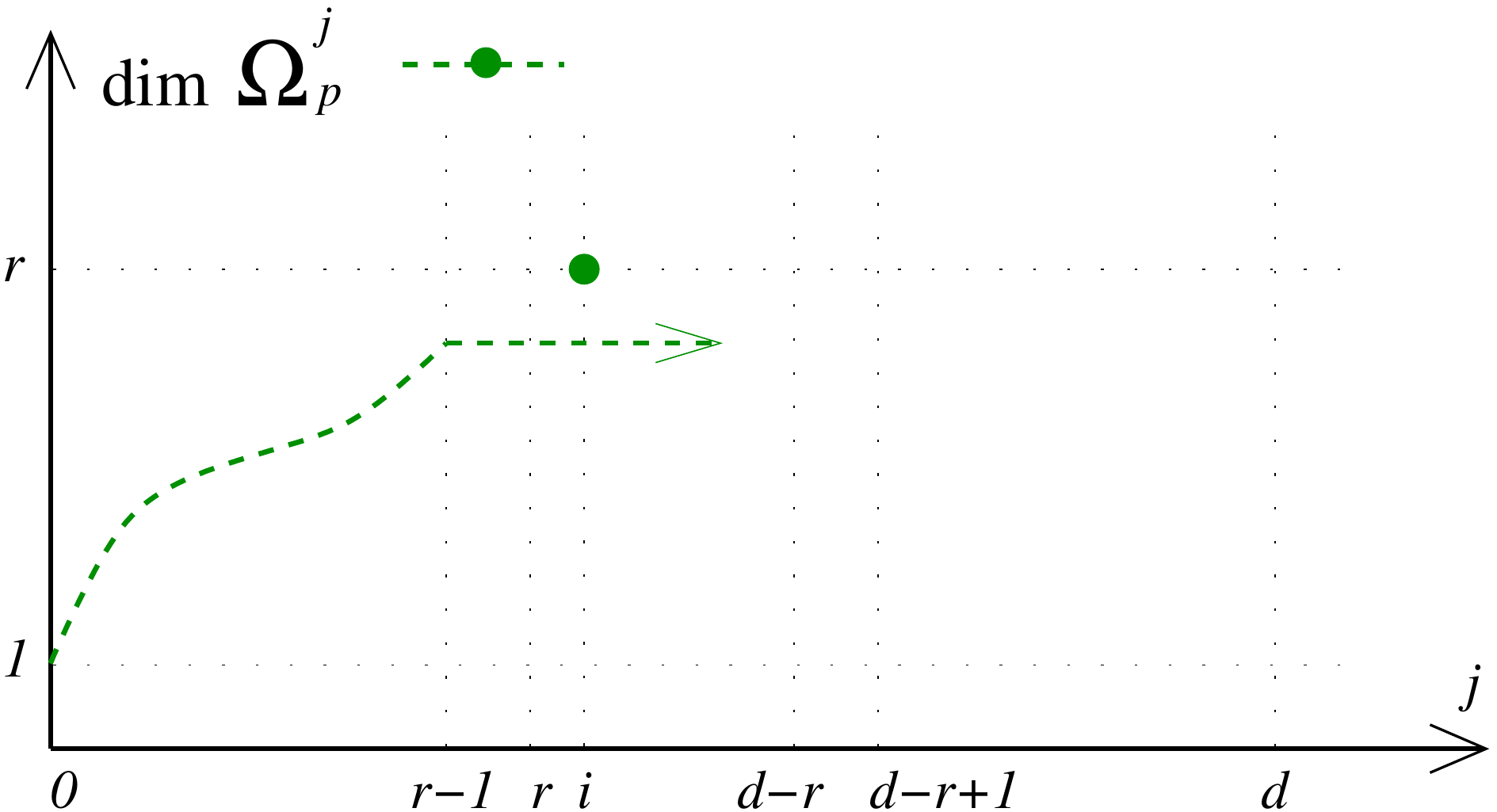}
\captionof{figure}{
If $\dim \Omega_p^{r-1} \le r-1$,
then, by Corollary~\ref{cor_Macaulays_bound_for_the_growth_of_Hilbert_function},
  we obtain $\dim \Omega_p^{\indexb} \le r-1$ for $\indexb \ge r-1$,
a contradiction with $\dim \Omega_p^{\indexa}=r$.
}\label{figure_with_dim_Omega_p_3}
\end{minipage}
&
\begin{minipage}[t]{0.47\textwidth}
\includegraphics[width=\textwidth]{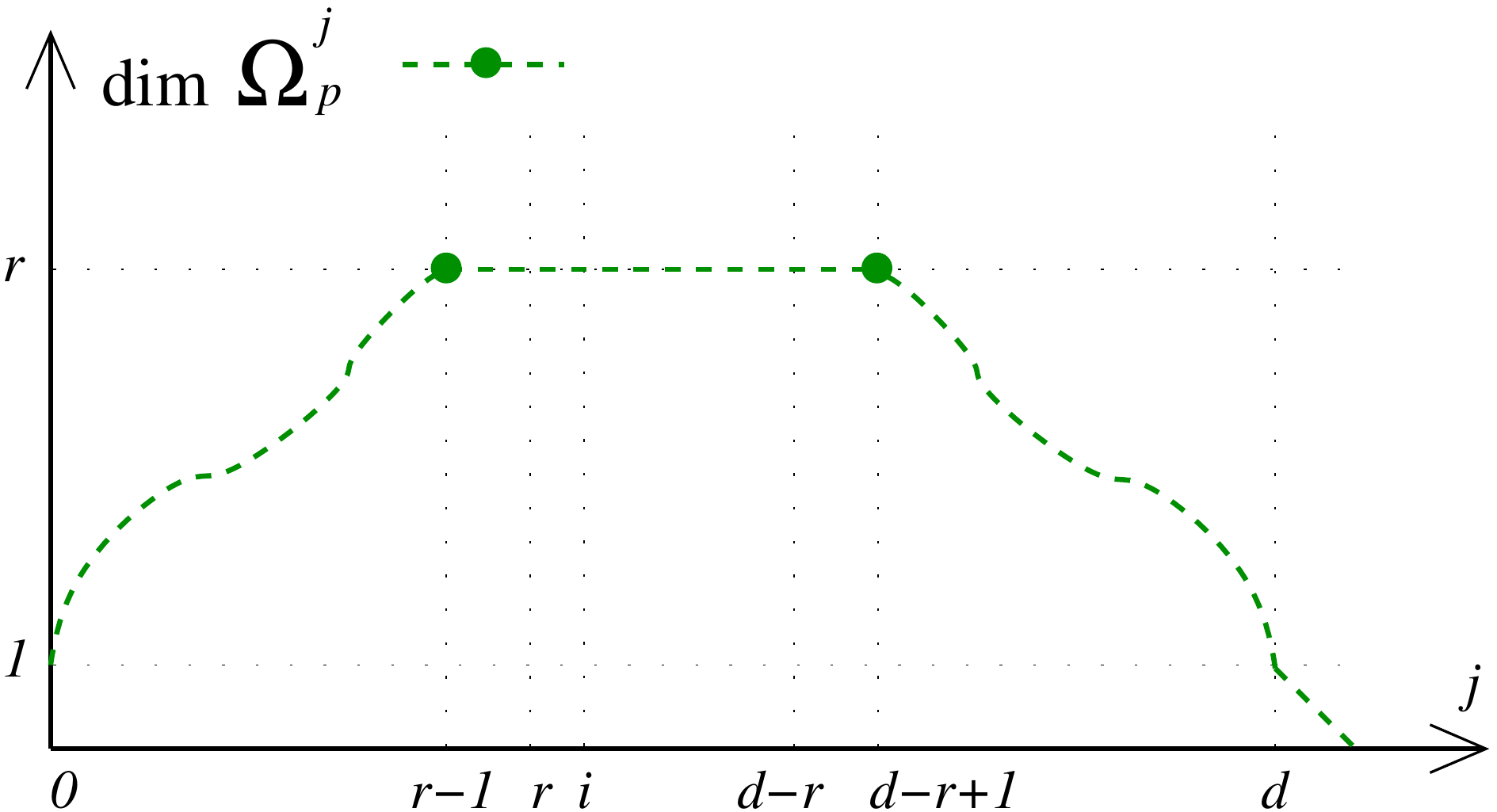}
\captionof{figure}{
Therefore $\dim \Omega_p^{r-1} =r$ and also $\dim \Omega_p^{d-r+1} =r$.
By the arguments under Figures~\ref{figure_with_dim_Omega_p_1} and \ref{figure_with_dim_Omega_p_2},
also $\dim \Omega_p^{\indexb} =r$ for $\indexb \in \setfromto{r-1}{d-r+1}$.
Thus the middle part is also proved.}\label{figure_with_dim_Omega_p_4}
\end{minipage}
\end{tabular}
\end{figure}

\begin{lemma}\label{lemma_hilbert_function_of_Omega_p_is_nice}
   Suppose $d \ge 2r$, $\indexa \in \setfromto{r}{d-r}$ and $\dim \Omega_p^{\indexa} = r$.
   Then:
   \begin{itemize}
     \item $\dim \Omega_p^{\indexb} \le \dim \Omega_p^{\indexb+1}$ for $\indexb \in \setfromto{0}{r-2}$;
     \item $\dim \Omega_p^{\indexb} = r$ for $\indexb \in \setfromto{r-1}{d-r+1}$;
     \item $\dim \Omega_p^{\indexb} \ge \dim \Omega_p^{\indexb+1}$ for $\indexb \in \setfromto{d-r+1}{d-1}$.
   \end{itemize}
   That is, the Hilbert function of $\Omega_p$ behaves as illustrated on Figure~\ref{figure_with_dim_Omega_p_4}.
\end{lemma}

\begin{prf}
By symmetry $\dim \Omega_p^{\indexb} = \Omega_p^{d-\indexb}$
  from Proposition~\ref{prop_basic_properties_of_Omega_p_and_I_p}\ref{item_Omega_is_perpendicular_to_I},
  we may suppose without loss of generality that $\indexa\le \frac{d}{2}$.
The argument follows in the captions of Figures~\ref{figure_with_dim_Omega_p_1}--\ref{figure_with_dim_Omega_p_4}.
\end{prf}

Suppose $\Psi$ is a graded algebra  generated in degree $1$ and
  $ \Psi = \Sym V^* / \ccJ$ for a homogeneous ideal $\ccJ$,
  which is generated in degrees at most $r$.
The Gotzmann's Persistence Theorem \cite{gotzmann_persistence_theorem} or \cite[Thm~3.8]{green_generic_initial_ideals}
states, that if the Macaulay bound on growth of Hilbert function of $\Psi$
   is attained at $r$, then it is attained for all $\indexb \ge r$.
The following is the special case we are going to use:

\begin{cor}\label{cor_Gotzmann_Persistence_Theorem}
  Suppose $\Psi$ is a graded algebra generated in degree $1$ and
    $ \Psi = \Sym V^* / \ccJ$ for a homogeneous ideal $\ccJ$,
    which is generated in degrees at most $r$ for some integer $r$.
  Suppose $\dim \Psi^{r} = \dim \Psi^{r+1} = r$.
  Then $\dim \Psi^{\indexb} = r$ for all $\indexb \ge r$.
\end{cor}

The next lemma is the main technical step in the proof of Theorem~\ref{thm_our_version_of_Eisenbuds_for_Pn}.

\begin{lemma}\label{lemma_J_is_saturated}
  Suppose $d\ge 2r$ and $r \le \indexa \le d-r$.
  Suppose $[p] \in \Ranklocus{r}{\indexa}{d-\indexa}{\PP V}$.
  Define $\ccJ$ to be the homogeneous ideal generated by the first $r$ degrees of $\annihp$,
    that is $\ccJ$ is generated by $\bigoplus_{\indexb=0}^{r} \annihp^{\indexb}$.
  Then $\ccJ$ is saturated,
    it defines a zero-dimensional scheme $R \subset \PP V$ of degree $\dim \Omega_p^{\indexa} \le r$,
    and the ideal $\annihp$ has no generators in degrees $\fromto{r+1}{d-r+1}$.
\end{lemma}

\begin{prf}
   Without loss of generality, assume $[p] \notin \Ranklocus{r-1}{\indexa}{d-\indexa}{\PP V}$,
      that is $\dim \Omega_p^{\indexa} = r$.
   Let $R \subset \PP V$ be the scheme defined by  $\ccJ$,
     and let $\Psi:= \Sym V^* / \ccJ$.
   We claim that for $\indexb \ge r$, we have $\dim\Psi^{\indexb} = r$.
   To obtain this we use Lemma~\ref{lemma_hilbert_function_of_Omega_p_is_nice},
      the Macaulay's bound and Gotzmann's Persistence Theorem.
   The details of the argument are illustrated on
   Figures~\ref{figure_with_dim_Psi_1}--\ref{figure_with_dim_Psi_4} and are following:
\begin{figure}[htb]
\centering
\textit{Figures~\ref{figure_with_dim_Psi_1}--\ref{figure_with_dim_Psi_4}.
        Illustration for the proof of Lemma~\ref{lemma_J_is_saturated}.}
\bigskip
\begin{tabular}{cc}
\begin{minipage}[t]{0.45\textwidth}
\includegraphics[width=\textwidth]{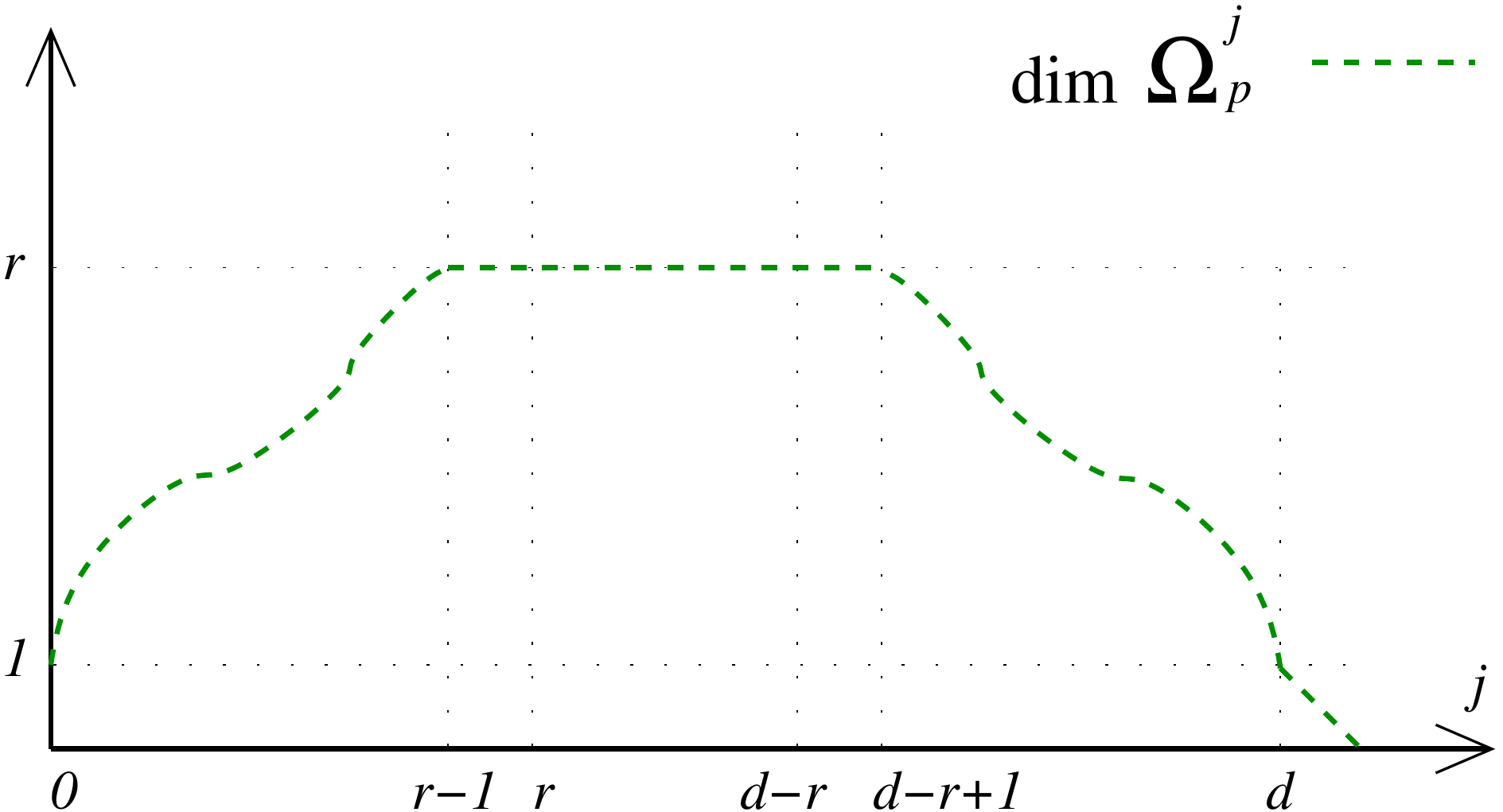}
\captionof{figure}{The graph of $\dim \Omega_p^{\indexb}$ obtained from Lemma~\ref{lemma_hilbert_function_of_Omega_p_is_nice}.}
\label{figure_with_dim_Psi_1}
\end{minipage}
&
\begin{minipage}[t]{0.45\textwidth}
\includegraphics[width=\textwidth]{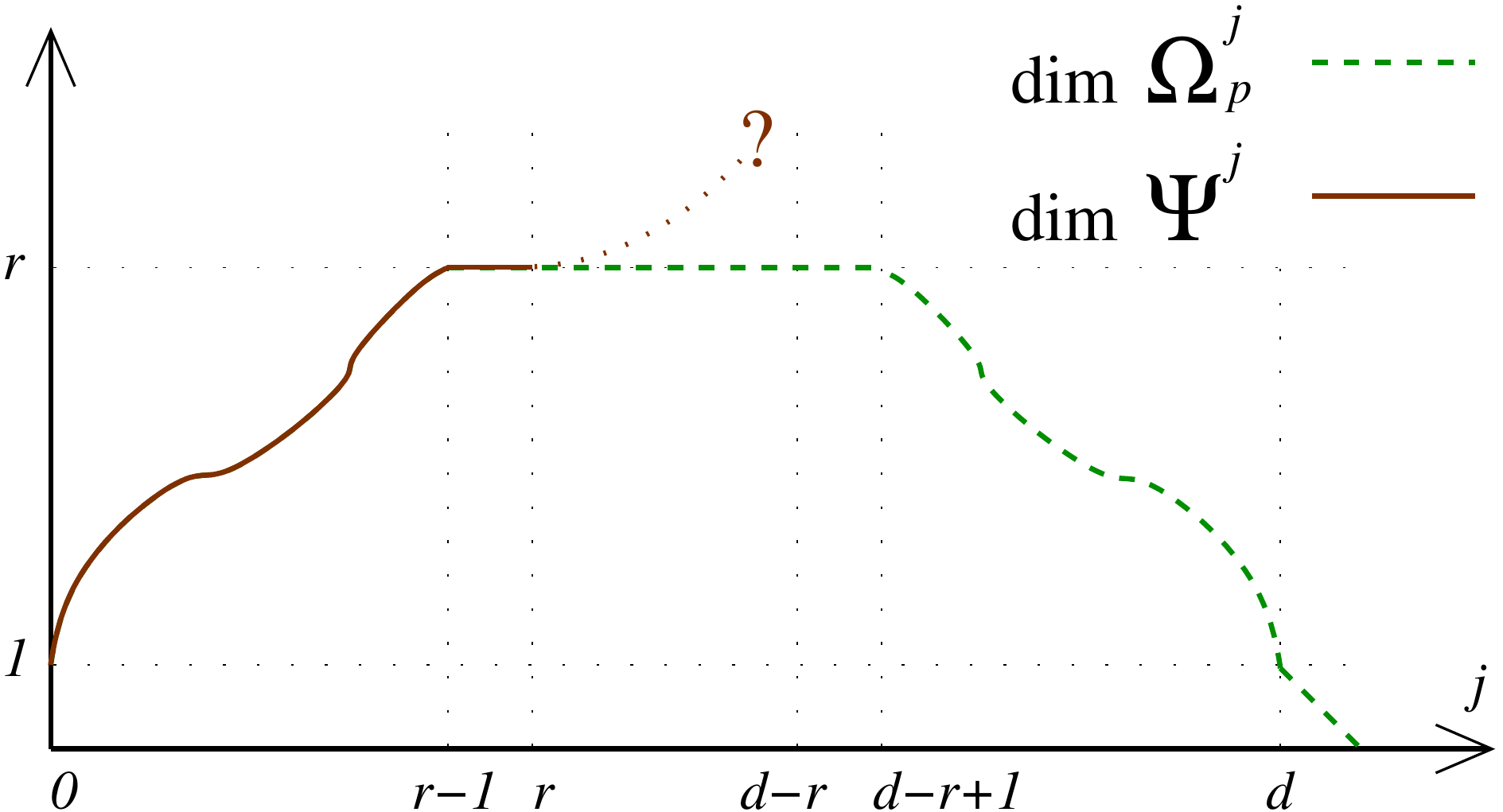}
\captionof{figure}{
By definition of $\ccJ$, Hilbert function of $\Psi = \Sym V^* /\ccJ$
  agrees with the  Hilbert function of $\Omega_p$ up to degree $r$
  --- see \ref{item_in_prf_Thm_our_Eisenbud_Hilb_functions_equal_at_the_begining}.
Moreover $\dim \Psi^{\indexb} \ge \dim \Omega_p^{\indexb}$ for all $\indexb$
  --- see \ref{item_in_prf_Thm_our_Eisenbud_Hilb_functions_ge_later}.}\label{figure_with_dim_Psi_2}
\end{minipage}
\\
\\
% \end{tabular}
% \end{figure}
% \begin{figure}[htb]
% \centering
% \begin{tabular}{cc}
\begin{minipage}[t]{0.45\textwidth}
\includegraphics[width=\textwidth]{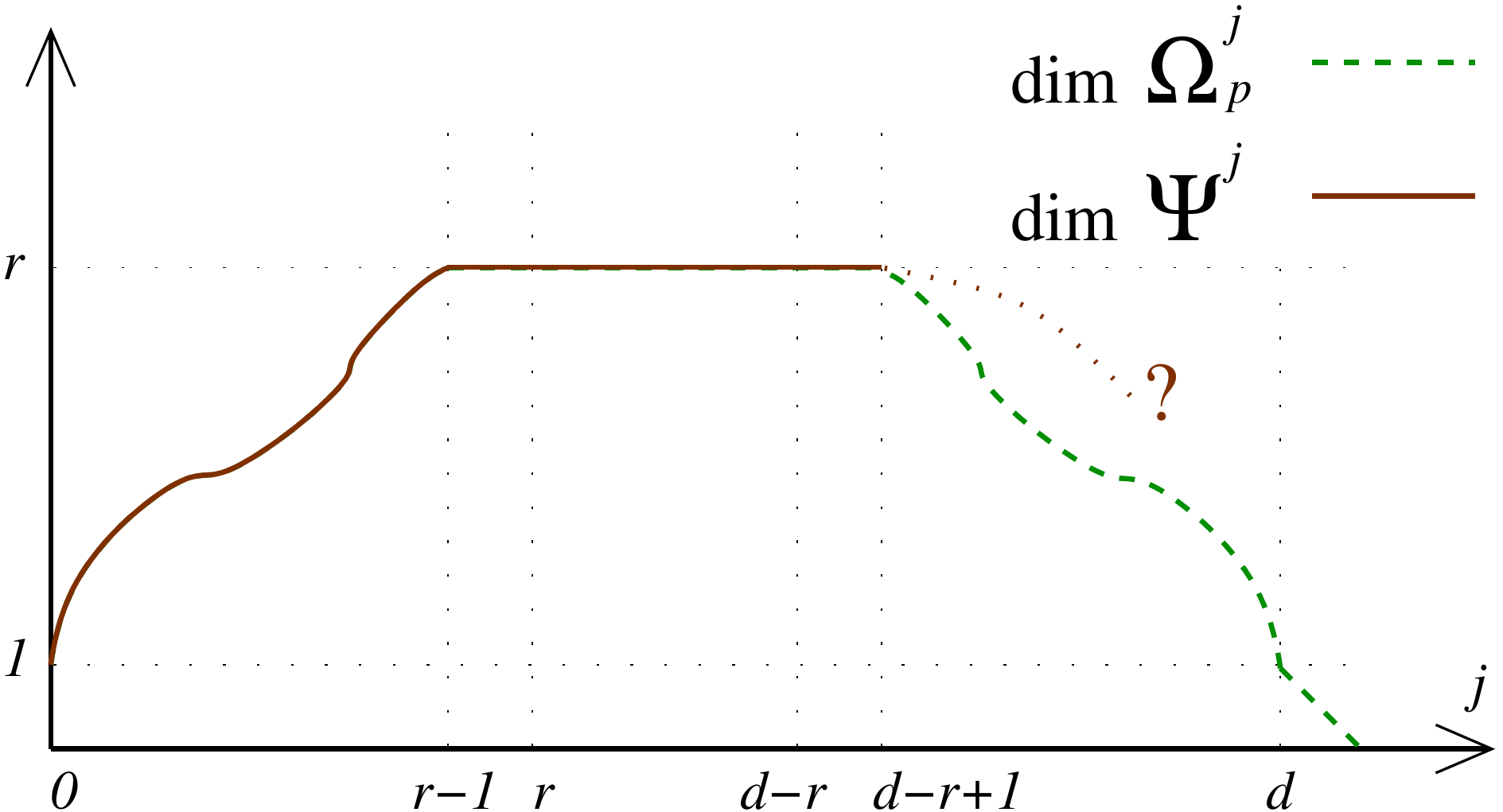}
\captionof{figure}{Since
$\dim \Psi^{r} = r$,
  by Corollary~\ref{cor_Macaulays_bound_for_the_growth_of_Hilbert_function},
  we obtain $\dim \Psi^{\indexb} \le r$ for $\indexb \ge r$
  --- see \ref{item_in_prf_Thm_our_Eisenbud_Hilb_functions_macaulay_bound}.}
  \label{figure_with_dim_Psi_3}
\end{minipage}
&
\begin{minipage}[t]{0.45\textwidth}
\includegraphics[width=\textwidth]{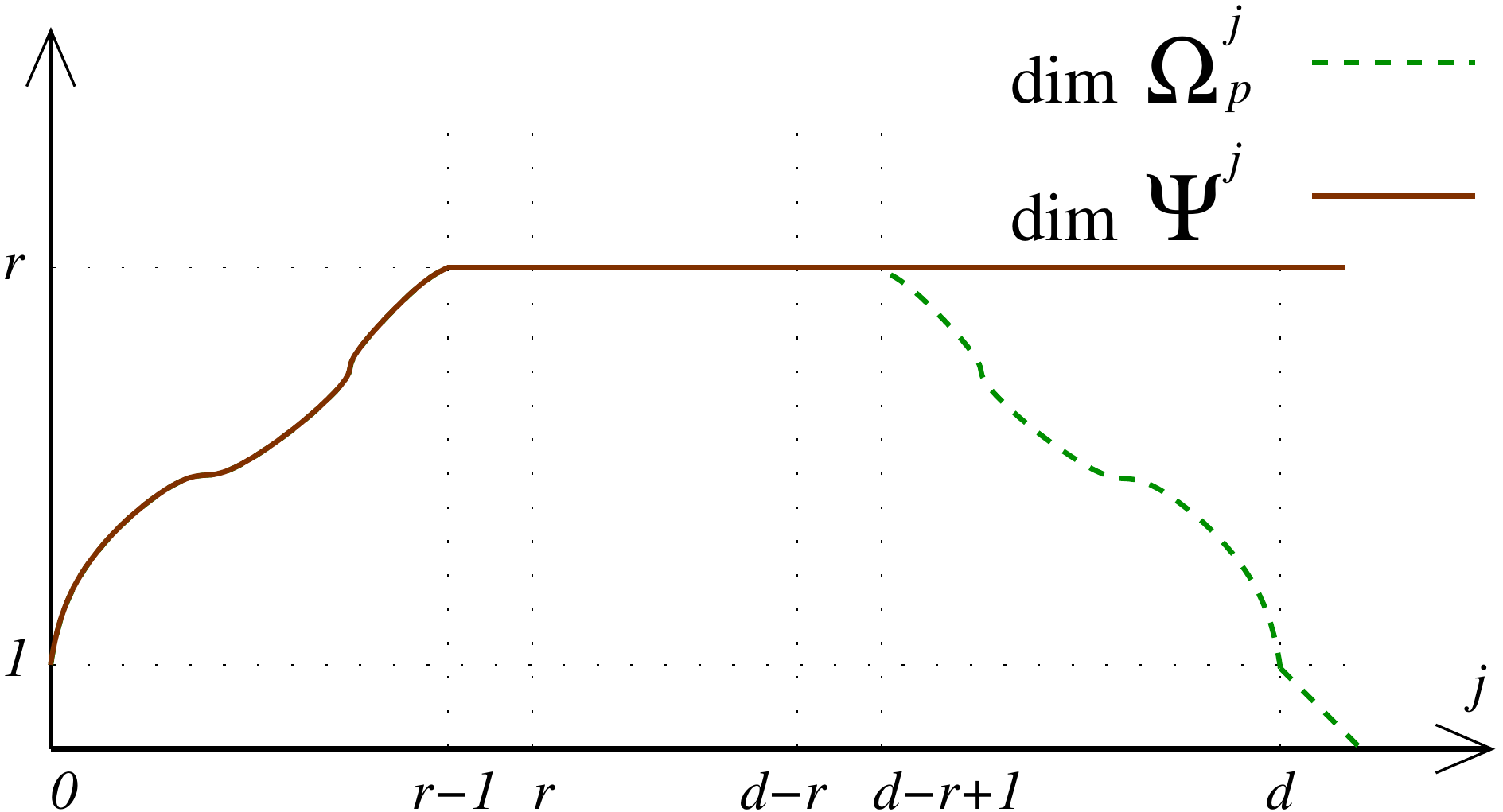}
\captionof{figure}{
By Corollary~\ref{cor_Gotzmann_Persistence_Theorem} the growth persists and
  $\dim \Psi^{\indexb} = r$ for $\indexb \ge r$ --- see \ref{item_in_prf_Thm_our_Eisenbud_Persistence}.}
  \label{figure_with_dim_Psi_4}
\end{minipage}
\end{tabular}
\end{figure}
  \renewcommand{\theenumi}{(\textbf{\textsc{\alph{enumi}}})}
   \begin{enumerate}
     \item  \label{item_in_prf_Thm_our_Eisenbud_Hilb_functions_equal_at_the_begining}
           $\dim \Psi^{\indexb} = \dim \Omega_p^{\indexb}$ for $\indexb \in \setfromto{1}{r}$, by definitions;
     \item  \label{item_in_prf_Thm_our_Eisenbud_Hilb_functions_ge_later}
           $\dim \Psi^{\indexb} \ge \dim \Omega_p^{\indexb}$ for all $\indexb$,
             because $\ccJ \subset \annihp$ --- see Figure~\ref{figure_with_dim_Psi_2};
     \item  \label{item_in_prf_Thm_our_Eisenbud_Omega_epsilon_ge_r}
           $\dim \Psi^{\indexb} \ge r$ for $\indexb \in \setfromto{r}{d-r+1}$,
             by \ref{item_in_prf_Thm_our_Eisenbud_Hilb_functions_ge_later}
               and Lemma~\ref{lemma_hilbert_function_of_Omega_p_is_nice};
     \item  \label{item_in_prf_Thm_our_Eisenbud_Omega_r_eq_r}
           $\dim \Psi^{r} = r$
              by \ref{item_in_prf_Thm_our_Eisenbud_Hilb_functions_equal_at_the_begining}
               and Lemma~\ref{lemma_hilbert_function_of_Omega_p_is_nice}.
     \item  \label{item_in_prf_Thm_our_Eisenbud_Hilb_functions_macaulay_bound}
           $\dim \Psi^{\indexb} \le r$ for $\indexb \ge r$
              by \ref{item_in_prf_Thm_our_Eisenbud_Omega_r_eq_r}
               and Corollary~\ref{cor_Macaulays_bound_for_the_growth_of_Hilbert_function}
               --- see Figure~\ref{figure_with_dim_Psi_3};
     \item  \label{item_in_prf_Thm_our_Eisenbud_Omega_r_plus_1_eq_r}
           $\dim \Psi^{r+1} = r$
              by \ref{item_in_prf_Thm_our_Eisenbud_Omega_epsilon_ge_r}
               and \ref{item_in_prf_Thm_our_Eisenbud_Hilb_functions_macaulay_bound};
               note that here we use that $d\ge 2r$, or equivalently $r+1 \le d-r+1$;
     \item \label{item_in_prf_Thm_our_Eisenbud_Persistence}
           $\dim \Psi^{\indexb} = r$  for all $\indexb \ge r$,
              by \ref{item_in_prf_Thm_our_Eisenbud_Omega_r_eq_r},
                  \ref{item_in_prf_Thm_our_Eisenbud_Omega_r_plus_1_eq_r}
               and Corollary~\ref{cor_Gotzmann_Persistence_Theorem}
	       --- see Figure~\ref{figure_with_dim_Psi_4}.
   \end{enumerate}
   \renewcommand{\theenumi}{(\roman{enumi})}
   Thus by \ref{item_in_prf_Thm_our_Eisenbud_Persistence}
     the Hilbert polynomial of $\Psi$ is equal to the constant polynomial $r$
     and $\Psi^{\indexb} = \Omega_p^{\indexb}$ for
   $\indexb \in \setfromto{r}{d-r+1}$.
   So $\dim R =0 $ and $\deg R =r$ and $\annihp$ has no generators in degrees $\setfromto{r+1}{d-r+1}$.
   Thus  the Lemma is proved,
      except for the statement that $\ccJ$ is saturated.

   Since $\ccJ$ agrees with $\annihp$ in the first $r$ degrees,
     by Proposition \ref{prop_basic_properties_of_Omega_p_and_I_p}\ref{item_I_p_is_delta_saturated}
     we have
     \[
	\ccJ^{\indexb}
	  = \set{\alpha \in S^{\indexb} V^* \mid S^{r-\indexb} V^* \cdot \alpha  \subset \ccJ^{r}}.
     \]
     for $\indexb \le r$.
   In the language of \cite[p.~298]{iarrobino_kanev_book_Gorenstein_algebras},
     this means that $\ccJ$ is \emph{the ancestor ideal} of $\ccJ^r$.
   From \cite[Cor.~C.18]{iarrobino_kanev_book_Gorenstein_algebras}
     we conclude that $\ccJ$ is saturated.
 \end{prf}

This enables us to prove Theorems~\ref{thm_our_version_of_Eisenbuds_for_Pn}
   and \ref{thm_our_version_of_Eisenbuds_for_Pn_effective_version}.
These theorems state that for sufficiently large $d$,
  the cactus variety $\cactus{r}{v_d(\PP V)}$
  is set-theoretically cut out by catalecticants and describe explicitly
  how to obtain the scheme $R$ such that $p \in \langle v_d (R) \rangle$.

\begin{prf}[ of Theorem~\ref{thm_our_version_of_Eisenbuds_for_Pn_effective_version}]
    Our assumptions are  $d\ge 2r$ and $r \le \indexa  \le d-r$
     and $p\in S^d V$ is such that $[p] \in \Ranklocus{r}{\indexa}{d-\indexa}{\PP V}$.
   $\ccJ$ is the homogeneous ideal generated by the first $\indexc$ degrees of $\annihp$,
     where $\indexc$ is any number such that $r \le \indexc \le d-r+1$ and
   $R \subset \PP V$ be the scheme defined by $\ccJ$.

   Lemma~\ref{lemma_J_is_saturated} implies that $\ccJ$ is saturated,
     does not depend on the choice of $\indexc$,
     $\dim R = 0$ and $\deg R \le r$.
   Thus part~\ref{item_J_is_saturated_and_independent_of_choice}
     is proved and it remains to prove  $[p] \in \langle v_d (R) \rangle$
     and the uniqueness of \ref{item_R_is_unique}.

   Since $\ccJ$ is saturated, $\langle v_d (R) \rangle = \PP\left( (\ccJ^d)^\perp\right)$.
   Since $\ccJ \subset \annihp$, we have $p \in (\annihp^d)^{\perp} \subset (\ccJ^d)^\perp$
     and thus $[p] \in \langle v_d (R) \rangle$.
   Thus \ref{item_R_has_the_right_dim_deg_and_has_p_in_span} is proved.

   To prove the minimality of \ref{item_R_is_unique},
     note that there is no zero-dimensional scheme $R' \subset \PP V$
     of degree strictly less than $\dim \Omega_p^{\indexa}$, such that $[p] \in \langle v_d(R') \rangle$
     by Proposition~\ref{prop_cactus_is_contained_in_rank_locus}.
   Thus the uniqueness follows from \cite[Cor.~2.2.1]{jabu_ginensky_landsberg_Eisenbuds_conjecture}.
\end{prf}

\begin{prf}[ of Theorem~\ref{thm_our_version_of_Eisenbuds_for_Pn}]
   We want to prove, that set-theoretically:
   \[
     \cactus{r}{v_d(\PP V)} = \Ranklocus{r}{\indexa}{d-\indexa}{\PP V}.
   \]
   The inclusion $\subset$ follows from Proposition~\ref{prop_cactus_is_contained_in_rank_locus}.
   The inclusion $\supset$ follows from Theorem~\ref{thm_our_version_of_Eisenbuds_for_Pn_effective_version}.
\end{prf}

\begin{prf}[ of Corollary~\ref{cor_decomposition_for_points_on_pure_secants}]
  If $p = {v_1}^d + \dotsb + {v_r}^d$ for some $v_i \in V$,
    then set $Q$ to be the reduced scheme $\setfromto{[v_1]}{[v_r]} \subset \PP V$.
  By Theorem~\ref{thm_our_version_of_Eisenbuds_for_Pn_effective_version}\ref{item_R_is_unique}
    we have $R \subset Q$.
  In particular, $R$ is reduced, since $Q$ is reduced.

  Now suppose $R$ is reduced, supported at $\setfromto{[v_1]}{[v_r]} \subset \PP V$
    (if $\deg R < r$, then take $v_{\deg R+1} = \dotsb = v_{r}=0 \in V$).
  Since $[p] \in \langle v_d (R) \rangle$,
    $p = a_1 {v_1}^d + \dotsb + a_r {v_r}^d$ for some $a_i \in \CC$.
  Rescaling the $v_i$ to absorb the constants, we obtain   $p =  {v_1}^d + \dotsb + {v_r}^d$.
\end{prf}

\section{Smoothability of Gorenstein schemes}\label{section_smoothability_of_Gorenstein}

The subject of smoothability of zero-dimensional Gorenstein schemes is intensively studied,
   see \cite[p.~221]{iarrobino_kanev_book_Gorenstein_algebras} or
   \cite{casnati_notari_irreducibility_Gorenstein_degree_9},
   \cite{casnati_notari_irreducibility_Gorenstein_degree_10},
   for  some of the recent results.

Suppose $X$ is a projective variety and $R \subset X$ is a zero-dimensional Gorenstein scheme of degree $r$.
If $X$ is singular and $R$ is supported in the singular locus,
  then it is relatively easy to give examples with  $R$ non-smoothable in $X$
  --- see for instance \cite[\S3]{jabu_ginensky_landsberg_Eisenbuds_conjecture},
  where such examples with $r=2$ and $X$ a singular curve are constructed.
If however $X$ is smooth, then  $R$ is smoothable in $X$ if and only if $R$ is smoothable,
see Proposition~\ref{prop_smothable_in_X_iff_smoothable_in_Y}.

The two propositions below summarise the known results about smoothability of Gorenstein zero-dimensional schemes.

\begin{prop}\label{prop_when_Gorenstein_are_smoothable}
  Let $R$ be a Gorenstein zero-dimensional scheme of length $r$ and embedding dimension $n$.
  If either $r \le 10$ or $n \le 3$, then $R$ is smoothable.
\end{prop}
For proof in the case $n\le 3$, see \cite[Cor.~2.6]{casnati_notari_irreducibility_Gorenstein_degree_9},
in the case $r\le 10$ see \cite{casnati_notari_irreducibility_Gorenstein_degree_10} and references therein.
See Section~\ref{section_improving_bounds} for a brief description of further work in progress in this direction.

\begin{prop}\label{prop_when_Gorenstein_are_not_smoothable}
  Suppose one of the following holds:
  \begin{itemize}
    \item  $n \ge 6, r \ge 14$ or
    \item  $n = 5, r \ge 42$ or
    \item  $n = 4, r \ge 140$.
  \end{itemize}
  Then there exists a zero-dimensional, degree $r$, non-smoothable Gorenstein scheme  $R\subset \CC^n$.
\end{prop}

\begin{prf}
  It is enough to construct the three extremal cases
  $(n,r) = (6, 14)$, $(5,42)$ and $(4,140)$,
  as we can increase $n$ by reembedding $R$ into a higher dimensional space,
  and we can increase $r$ by adding to $R$ disjoint points.
  We claim that for these $(n,r)$  there exists
  a zero-dimensional local Gorenstein algebra $\Psi$,
  such that $\Spec \Psi$ is not smoothable.
  We construct $\Psi = \Sym V^*/\ccI$
    by setting $\ccI:=\annihp$ for a general, not necessarily homogeneous polynomial of degree $\indexb$ in $n$ variables.
  Here $\indexb = 3$ for $n=6$, $\indexb=5$ for $n=5$, and $\indexb=9$ for $n=4$.

  Iarrobino observed that the case $(n,\indexb) = (6, 3)$ gives rise to non-smooth\-able scheme of degree $14$
    by calculating the tangent space to the Hilbert scheme at the point represented by such scheme
    (the tangent space is too small) ---
    see \cite[Lem.~6.21]{iarrobino_kanev_book_Gorenstein_algebras}
    or \cite[\S4, p11]{casnati_notari_irreducibility_Gorenstein_degree_10}.

  Since $p$ is general, by \cite[Thm~3.31]{emsalem_iarrobino_small_tangent_space}
    and \cite[Thm~1D]{iarrobino_compressed_artin},
    the algebra $\Psi$ is \emph{compressed}, that is, it has the maximal possible length:
   \[
     r := \dim_{\CC} \Psi = 2\sum_{i=0}^{\lfloor \frac{\indexb}{2} \rfloor} \dim S^{i} \CC^n
         = 2 \binom{\lfloor \frac{\indexb}{2} \rfloor + n}{n}.
   \]
   Here we only wrote the formula for odd $\indexb$.
   In particular, in our cases $(n,\indexb) = (5, 5)$ and  $(4,9)$  we get,  respectively,
     $r = 42$ and  $140$.

   Let $Z(n,\indexb)$ be the parameter space for all Gorenstein subschemes in $\CC^n$
   supported at $0 \in \CC^n$ with socle degree $\indexb$ and of maximal length $r$.
   Its dimension is calculated in \cite[Thm~2]{iarrobino_compressed_artin}:
   \[
     \dim Z(n, \indexb) = \binom{n+\indexb}{\indexb} - r.
   \]
   Thus in the cases $(n,\indexb, r) = (5, 5, 42)$ and  $(4,9, 140)$
     we have $\dim Z(n, \indexb) = 210$ and $575$ respectively.

   Let $\operatorname{H}_r(\CC^n)$ be the smoothable component,
     that is the irreducible scheme para\-metri\-sing smoothable subschemes of $\CC^n$ of length $r$.
   Its dimension is $nr$ and a general point represents general $r$ distinct points.
   For $(n,r) = (5,42)$ and $(4,140)$ we obtain, respectively,
     $\dim \operatorname{H}_r(\CC^n) = 210$ and $560$.
   Thus in these two cases $\dim Z(n, \indexb) \ge \dim \operatorname{H}_r(\CC^n)$
     and therefore it is impossible that $Z(n, \indexb) \subset \operatorname{H}_r(\CC^n)$ ---
     it is obvious, when $>$ holds, and if $=$ holds,
     then, since $\operatorname{H}_r(\CC^n)$ is irreducible,
     it could only happen if $Z(n, \indexb) \subset \operatorname{H}_r(\CC^n)$ is dense.
   But a general point of $Z(n, \indexb)$ represents a scheme supported at $0$,
     not at $r$ distinct points.
   Thus a general scheme in $Z(n, \indexb)$ is not smoothable,
     and $ R:= \Spec \Psi$ is not smoothable.
\end{prf}

\begin{rem}
  We underline that the role of $\indexb$ in the proof of Proposition~\ref{prop_when_Gorenstein_are_not_smoothable}
    is different than that of $d$
    in Sections~\ref{section_flattenings_and_Gorenstein_Artin_algebras}--\ref{section_bounds_on_Hilb_functions}.
  As one example, in Theorems~\ref{thm_secants_are_cut_out_by_catalecticants} and \ref{thm_counterexamples_to_Eisenbud}
    we assume $d\ge 2r-1$ or $d\ge 2r$, whereas here $j$ is very small compared to $r$.
  In particular, in the proof of Theorem~\ref{thm_counterexamples_to_Eisenbud} below
    we are going to \emph{reembed} the scheme $R$ constructed above.
\end{rem}

\begin{rem}
   To construct the non-smoothable example of Iarrobino in the case $(n,r) = (6, 14)$
     it is enough to use as $p$ a general \emph{homogeneous} polynomial of degree $\indexb=3$.
   One reason for that is the following: for a general non-homogeneous $p$ of degree $3$ the algebra $\Omega_p$
     is compressed and therefore, by using the normal form of \cite[Thm~4.1]{elias_rossi_short_Gorenstein},
     it is isomorphic to a graded algebra $\Omega_{p'}$, for homogeneous $p'$.
   However, it is not known to the authors if to obtain the non-smoothable examples with $(n,r) = (5, 42)$ and $(4,140)$
     it is enough to consider general homogeneous $p$ of degree $\indexb = 5$ or $9$, respectively.
   The inequalities in our proof (coming from the dimension count)
     do not work if we restrict our attention to the homogeneous polynomials,
     and (although quite unlikely) it is potentially possible, that all homogeneous algebras $\Omega_p$
     are in the intersection of two components of the Hilbert scheme and therefore are smoothable.
\end{rem}

Now we can conclude the Theorems~\ref{thm_secants_are_cut_out_by_catalecticants} and \ref{thm_counterexamples_to_Eisenbud}
hold.

Theorem \ref{thm_secants_are_cut_out_by_catalecticants} says that in the range of $n$, $r$
  as in Proposition~\ref{prop_when_Gorenstein_are_smoothable}
the secant variety to the Veronese variety of sufficiently high degree is set-theoretically cut out by catalecticant minors.

\begin{prf}[ of Theorem~\ref{thm_secants_are_cut_out_by_catalecticants}]
Since $r\le 10$ or $n \le 3$,
  all zero-dimensional subschemes of $\PP V \simeq \PP^n$ of degree $\le r$ are smoothable,
  by Proposition~\ref{prop_when_Gorenstein_are_smoothable}.
Thus Theorem~\ref{thm_when_secant_equals_cactus}\ref{item_if_star_holds_then_sigma_eq_cactus}
   applies and $\sigma_r(v_d(\PP V)) = \cactus{r}{v_d(\PP V)}$.
Further, since  $d\ge 2r$ and $r \le \indexa \le d-r$,
  the assumptions in Theorem~\ref{thm_our_version_of_Eisenbuds_for_Pn}
  are satisfied and therefore
  $\cactus{r}{v_d(\PP V)} = \Ranklocus{r}{\indexa}{d-\indexa}{\PP V}$ as sets.
\end{prf}

On the contrary, Theorem~\ref{thm_counterexamples_to_Eisenbud}
  states that in the range of $n,r$ as in Proposition~\ref{prop_when_Gorenstein_are_not_smoothable},
  the secant variety to the Veronese variety of sufficiently high degree is not defined by the catalecticant minors.

\begin{prf}[ of Theorem~\ref{thm_counterexamples_to_Eisenbud}]
  By Proposition~\ref{prop_basic_properties_of_cactus}:
  \[
     \sigma_r(v_d(\PP V)) \subset \cactus{r}{v_d(\PP V)}
  \]
  and by Proposition~\ref{prop_cactus_is_contained_in_rank_locus}:
  \[
     \cactus{r}{v_d(\PP V)} \subset \bigcap_{\indexa=1}^{d-1} \Ranklocus{r}{\indexa}{d-\indexa}{\PP V}.
  \]

By Proposition~\ref{prop_when_Gorenstein_are_not_smoothable},
  there are non-smoothable Gorenstein zero-dimen\-sional subschemes of $\PP V$
  of degree $r$.
Thus by Theorem~\ref{thm_when_secant_equals_cactus}\ref{item_if_sigma_eq_cactus_then_star_holds},
  we must have $\sigma_r(v_d(\PP V)) \varsubsetneq \cactus{r}{v_d(\PP V)}$.
Therefore $ \sigma_r(v_d(\PP V)) \varsubsetneq \bigcap_{\indexa=1}^{d-1} \Ranklocus{r}{\indexa}{d-\indexa}{\PP V}$ as sets.
\end{prf}

\begin{example}
   Let $(n,\indexb,r) = (6,3, 14)$ or $(5,5, 42)$ or $(4,9, 140)$, and
    $V \simeq \CC^{n+1}$
   and let $f \in S^{\indexb} V$ be a general element.
   Also let $z \in V$ be any non-zero element.
   Suppose $d \ge 2r-1$ and set $p:=f \cdot z^{d-\indexb} \in S^d V$.
   Then $p$ is a polynomial,
     whose catalecticant matrices have all rank at most $r$,
     but $[p]$ is not on the $r$-th secant variety of $v_d(\PP V)$.
   This is because
      the polynomial $p$ was chosen in such a way that $[p] \in \langle v_d(R) \rangle$,
      for a Gorenstein scheme $R$ constructed in the proof of Proposition~\ref{prop_when_Gorenstein_are_not_smoothable},
   supported at $[z]$.
   Thus $[p] \in \cactus{r}{v_d(\PP V)} \subset \Ranklocus{r}{\indexa}{d-\indexa}{\PP V}$ for any $\indexa$.
   Also $p \notin \langle v_d(R') \rangle$ for any $R' \varsubsetneq R$ and since $R$ is not smoothable,
      by \cite[Cor.~2.2.1]{jabu_ginensky_landsberg_Eisenbuds_conjecture} we must have
      $[p] \notin \sigma_r(v_d(\PP V))$.
\end{example}

\section{Secants to Veronese reembeddings of a smooth variety}\label{section_arbitrary_X}

The aim of this section is to apply the theory developed in previous sections to arbitrary $X \subset \PP V$
  and strengthen the results obtained in \cite{jabu_ginensky_landsberg_Eisenbuds_conjecture}.
These results are motivated by a question of Eisenbud
formulated in several different versions, see~\cite[\S1.2]{jabu_ginensky_landsberg_Eisenbuds_conjecture}.
In general, the problem is to determine when is the ideal of secant variety to a variety $X$ embedded by
  a ``sufficiently ample'' linear system defined by minors of a matrix with linear entries.
In the case of curves, a conjecture was formulated
in~\cite{eisenbud_koh_stillman_curves_of_high_degree}
  (see also \cite[Conj.~1.2.1]{jabu_ginensky_landsberg_Eisenbuds_conjecture}).
Later, the question in the form as
  in~\cite[Question~1.2.2]{jabu_ginensky_landsberg_Eisenbuds_conjecture}
  was informally suggested by David Eisenbud.
In these two forms, the question is more explicit about what linear system one needs to consider
  and which matrix should give the equations.
It has been also modified by Jessica Sidman and Greg Smith in \cite[Conj.~1.2]{sidman_smith_determinantal}
(see also \cite[Conj.~1.2]{raicu_3_times_3_minors})
  --- they claim that there is even more linear systems giving the embedding,
  but they have not specified the matrix providing the equations.
Finally, a restricted version was formulated as a question by Joseph Landsberg, Adam Ginensky, and the second author
  in~\cite[Question~1.2.3]{jabu_ginensky_landsberg_Eisenbuds_conjecture},
  see also Question~\ref{conj_restricted_Eisenbud} below.
All these questions however are proved to have negative answers, if we allow $X$ to have sufficiently bad singularities,
  see~\cite[Thm~1.1.4]{jabu_ginensky_landsberg_Eisenbuds_conjecture}.
For smooth $X$, there is a significant number of positive results,
  see~\cite[\S1.2]{jabu_ginensky_landsberg_Eisenbuds_conjecture}.
Yet our results in this paper, provide examples when the answer to
  \cite[Question~1.2.2]{jabu_ginensky_landsberg_Eisenbuds_conjecture} is negative.
Namely the projective space or any smooth variety of dimension at least $4$ is such an example.
We also suspect, that for sufficiently high dimension of $X$,
  also \cite[Conj.~1.2]{sidman_smith_determinantal} is too strong.
If it was true, then this would give a determinantal criterion to test
  if a Gorenstein scheme is smoothable or not.
We suspect that such an easy criterion should not exist.

Thus we concentrate on the weakest known version of the question,
  which is \cite[Question~1.2.3]{jabu_ginensky_landsberg_Eisenbuds_conjecture}
  with added smoothness assumption:
\begin{question}\label{conj_restricted_Eisenbud}
  Let $X\subset \PP V$ be a smooth,
  irreducible variety, and fix $r\in \NN$.
  Does there exists infinitely many $d$ such that
    $\sigma_r(v_d(X)) $ is cut out set-theoretically by equations of  degree $r+1$?
\end{question}

We are now ready to prove Theorem~\ref{thm_our_version_of_Eisenbuds_for_X}
  and Corollary \ref{cor_Eisenbuds_holds_for_X}.
The latter,
  by Proposition~\ref{prop_when_Gorenstein_are_smoothable},
  provides the affirmative answer to Question~\ref{conj_restricted_Eisenbud}
  in cases $\dim X \le 3$ or $r \le 10$.

\begin{prf}[ of Theorem~\ref{thm_our_version_of_Eisenbuds_for_X}]
  We suppose $X \subset \PP V$, $d \ge 2r$, $r \le \indexa \le d-r$
     and $d \ge Got(h_X) + r -1$.
  We claim the following equality of sets holds:
  \[
    \cactus{r}{v_d(X)} = \Ranklocus{r}{\indexa}{d-\indexa}{X} :=
                     \Ranklocus{r}{\indexa}{d-\indexa}{\PP V} \cap \langle v_d(X) \rangle
  \]

  The inclusion $\cactus{r}{v_d(X)} \subset \Ranklocus{r}{\indexa}{d-\indexa}{X}$ always holds.
  This is because
    $\cactus{r}{v_d(X)} \subset \langle v_d(X) \rangle$
    by Proposition~\ref{prop_basic_properties_of_cactus}\ref{item_sigma_subset_cactus} and
  \[
    \cactus{r}{v_d(X)} \subset \cactus{r}{v_d(\PP V)} \subset \Ranklocus{r}{\indexa}{d-\indexa}{\PP V}
  \]
    (the first inclusion holds by Proposition~\ref{prop_basic_properties_of_cactus}\ref{item_sigma_subset_cactus},
       and the second by Proposition~\ref{prop_cactus_is_contained_in_rank_locus}).

  To prove $\cactus{r}{v_d(X)} \supset \Ranklocus{r}{\indexa}{d-\indexa}{X}$,
    suppose $p\in \Ranklocus{r}{\indexa}{d-\indexa}{X} = \Ranklocus{r}{\indexa}{d-\indexa}{\PP V} \cap \langle v_d (X) \rangle$.
  Due to Theorem~\ref{thm_our_version_of_Eisenbuds_for_Pn},
    our assumptions on $r$, $d$ and $\indexa$
    imply that $\Ranklocus{r}{\indexa}{d-\indexa}{\PP V} = \cactus{r}{v_d(\PP V)}$.
  Thus, by Proposition~\ref{properties_of_cactus_of_Veronese_reembeddings},
    there exists a zero-dimensional scheme $R \subset \PP V$
    of degree at most $r$ such that $p \in \langle v_d(R) \rangle$.
  That is, $p \in \langle v_d(R) \rangle \cap \langle v_d (X) \rangle$.
  Since $d\ge Got(h_x) +r -1$, by \cite[Lem.~1.1.2]{jabu_ginensky_landsberg_Eisenbuds_conjecture},
    one has $p \in \langle v_d(R \cap X) \rangle$.
  Since $v_d(R \cap X)$ is a zero-dimensional scheme of degree at most $\deg R \le r$,
    $\langle v_d(R \cap X) \rangle \subset \cactus{r}{v_d(X)}$
    by the definition of cactus variety \eqref{equ_define_cactus2}.
  Therefore $p \in \cactus{r}{v_d(X)}$
    and the inclusion $\cactus{r}{v_d(X)} \supset \Ranklocus{r}{\indexa}{d-\indexa}{X}$ is proved.
\end{prf}

\begin{prf}[ of Corollary~\ref{cor_Eisenbuds_holds_for_X}]
  If condition \ref{item_condition_on_smoothability_of_all_Gorenstein} holds,
    then by Theorem~\ref{thm_when_secant_equals_cactus}
    we have $\sigma_r(v_d(X)) = \cactus{r}{v_d(X)}$.
  Further, by Theorem~\ref{thm_our_version_of_Eisenbuds_for_X},
    we also have
    $\cactus{r}{v_d(X)} = \Ranklocus{r}{\indexa}{d-\indexa}{X}$ as sets.

  If \ref{item_condition_on_smoothability_of_all_Gorenstein} fails to hold,
    then, by Theorem~\ref{thm_when_secant_equals_cactus}
    and Proposition~\ref{prop_cactus_is_contained_in_rank_locus},
    we have
   \[
     \sigma_r(v_d(X)) \varsubsetneq \cactus{r}{v_d(X)}
     \subset \bigcap_{\indexa=1}^{d-1} \Ranklocus{r}{\indexa}{d-\indexa}{X}.
   \]
\end{prf}

\section{Improving the bounds}\label{section_improving_bounds}

In this section we briefly discuss to what extend the bounds in the theorems presented in Section~\ref{section_intro}
 are effective.

\subsection{Bounds on rank $r$ and dimension $n$}\label{section_bounds_on_r_and_n}

According to Gianfranco Casnati and Roberto Notari,
  every Gorenstein zero-dimensional scheme of degree $r=11$
  is smoothable and notes explaining this \cite{casnati_notari_irreducibility_Gorenstein_degree_11}
  should be available shortly.
Thus the bound $r \le 10$ in Proposition~\ref{prop_when_Gorenstein_are_smoothable}
  and Theorem~\ref{thm_secants_are_cut_out_by_catalecticants} can be replaced by $r \le 11$.

Cases $r=12$, and $r=13$ are also investigated, but it is too early to give any definite answer.
One of the most essential problems to resolve these two cases,
   is to determine if a general
   graded Gorenstein Artin algebra with Hilbert function $(1,5,5,1,0,0,\dots)$ is smoothable.
By theorems in Section~\ref{section_intro} for $\PP^5 = \PP V$ where $V$ has basis
   $\fromto{x_1}{x_5}, z$,
   this problem is equivalent to determining if $[f\cdot z^{d-3}] \in \sigma_{12}(v_d(\PP  V))$,
   where $f$ is a general homogeneous cubic in $\fromto{x_1}{x_5}$ and $d$ is sufficiently large (at least $24$).

As explained in Proposition~\ref{prop_when_Gorenstein_are_not_smoothable},
   there are known examples of non-smooth\-able $R \subset \CC^n$ of degree $r=14$, provided $n\ge 6$.
It is not known, if there exists a low degree example in $\CC^4$ or $\CC^5$,
   but we expect, that the bounds $r\ge 140$ for $n=4$ and $r \ge 42$ for $n=5$
   in Proposition~\ref{prop_when_Gorenstein_are_not_smoothable}
   (and thus in Theorem~\ref{thm_counterexamples_to_Eisenbud})
   are far not effective.

\subsection{Bounds on degrees $d$ and $\indexa$}

We start this subsection with an easy example.

\begin{example}
   Let $d\ge 2r$ and $\indexa < r$.
   Let $p = {v_1}^d +\dotsb + {v_{r+1}}^d \in S^d \CC^2$,
      where $v_i \in \CC^2$ are general points.
   Then $[p] \in \Ranklocus{r}{\indexa}{d-\indexa}{v_d(\PP^1)}$, but $[p] \notin \cactus{r}{v_d(\PP^1)}$.
\end{example}
Thus the bound on $\indexa$ in Theorems~\ref{thm_secants_are_cut_out_by_catalecticants},
   \ref{thm_our_version_of_Eisenbuds_for_Pn}, \ref{thm_our_version_of_Eisenbuds_for_Pn_effective_version},
   \ref{thm_our_version_of_Eisenbuds_for_X} is optimal, provided $d \ge 2r$.
However, if $d\le 2r-1$, then sometimes we still have the equality of sets
   $\cactus{r}{v_d(\PP V)} =\Ranklocus{r}{\indexa}{d-\indexa}{v_d(\PP V)} $
  for $\indexa = \lfloor \frac{d}{2} \rfloor$.
Below we sketch the proof, in the case $d=2r-1$, which is divided into two cases,
  one of which is essentially identical to the proof for $d\ge 2r$.

\begin{prop}
   Suppose $d=2r-1$, $\indexa = r-1$.
   Then $\cactus{r}{v_{2r-1}(\PP V)} =\Ranklocus{r}{r-1}{r}{v_{2r-1}(\PP V)}$ as sets.
   Thus, if in addition $n\le 3$ or $r \le 10$
      (or $r \le 11$ taking in account \S\ref{section_bounds_on_r_and_n}),
      then $\sigma_r(v_{2r-1}(\PP V)) =\Ranklocus{r}{r-1}{r}{v_{2r-1}(\PP V)}$ as sets.
\end{prop}

\begin{prf}[ (sketch)]
   The strategy of the proof is identical to the proof of Theorem~\ref{thm_our_version_of_Eisenbuds_for_Pn},
     except we need some replacements for
     Lemmas~\ref{lemma_hilbert_function_of_Omega_p_is_nice} and \ref{lemma_J_is_saturated}.

   Let $[p]\in \Ranklocus{r}{r-1}{r}{v_{2r-1}(\PP V)}$.
   If $[p]\in \Ranklocus{r-1}{r-1}{r}{v_{2r-1}(\PP V)}$,
     then we can use Theorem~\ref{thm_our_version_of_Eisenbuds_for_Pn},
     so suppose $\dim  \Omega_p^{r-1} =r$.
   An analogue of   Lemma~\ref{lemma_hilbert_function_of_Omega_p_is_nice} (with an identical proof) gives
     the unimodality of the Hilbert function of $\Omega_p$:
   \begin{itemize}
     \item $\dim \Omega_p^{\indexb} \le \dim \Omega_p^{\indexb+1}$ for $\indexb \in \setfromto{0}{r-2}$;
     \item $\dim \Omega_p^{r-1} = \dim \Omega_p^{r} = r$.
     \item $\dim \Omega_p^{\indexb} \ge \dim \Omega_p^{\indexb+1}$ for $\indexb \in \setfromto{r}{2r-1}$;
   \end{itemize}

   Note that we cannot have $\dim \Omega_p^{r-2} \le r-2 $,
     as then by Corollary~\ref{cor_Macaulays_bound_for_the_growth_of_Hilbert_function}
     we would have $\dim \Omega_p^{r-1} \le r-2$ contrary to the above.
   Thus we consider two cases:  $\dim \Omega_p^{r-2} = r $ or   $\dim \Omega_p^{r-2} = r-1$.

   If $\dim \Omega_p^{r-2} = r $, then by symmetry $\dim \Omega_p^{r+1} = r$,
     and exactly the same argument as in Lemma~\ref{lemma_J_is_saturated} works
     (because the conclusion of \ref{item_in_prf_Thm_our_Eisenbud_Omega_r_plus_1_eq_r} holds).

   Now suppose $\dim \Omega_p^{r-2} = r-1$.
   We claim this is only possible if $p \in S^{2r-1} W$ for some linear subspace $W \subset V$
     with $\dim W =2$.
   Let $\ccJ \subset \Sym V^*$ be the ideal generated by first $r-2$ gradings of $\annihp$
     and set $\Psi:=\Sym V^*/\ccJ$.
   Then the Macaulay's bound
     \cite[Thm~2.2(i), (iii)]{stanley_hilbert_functions_of_graded_algebras} or
     \cite[Thm~3.3]{green_generic_initial_ideals}
     gives $\dim \Psi^{r-1} \le r$.
   Since $\ccJ \subset \annihp$, we have   $\dim \Psi^{r-1} = r$.
   Thus the maximal growth occurs and we can apply the Gotzmann's Persistence Theorem
     \cite{gotzmann_persistence_theorem} or \cite[Thm~3.8]{green_generic_initial_ideals}.
   We have  $\dim \Psi^{\indexb} = \indexb+1$ for all $\indexb \ge r-2$.
   If $X\subset \PP V$ is the scheme defined by $\ccJ$,
     then the Hilbert polynomial of $X$ is $\indexb+1$.
   The only such scheme is $\PP^1 = \PP W$  for linear subspace $W\subset V$.
   The saturation $\ccJ^{sat}$ of $\ccJ$ is the ideal defining $X$.
   So the Hilbert function of $\Sym V^* / \ccJ^{sat}$ is equal identically to $\indexb +1$,
     and thus it agrees with  $\dim \Psi^{\indexb}$ for all $\indexb \ge r-2$.
   In particular, $(\ccJ^{sat})^{2r-1}  = \ccJ^{2r-1} \subset \annihp^{2r-1}$
     and thus $[p] \in \PP((\ccJ^{sat})^{2r-1})^{\perp} = \langle v_{2r-1} (X) \rangle = \PP (S^{2r-1} W)$
     as claimed.

   But $\sigma_r (v_{2r-1} (\PP^1)) = \langle v_{2r-1} (\PP^1 )\rangle$,
     so $[p] \in \sigma_r(v_{2r-1} (\PP V)) \subset \cactus{r}{v_{2r-1} (\PP V)}$.
\end{prf}

In the proof above, we use the Macaulay's bound and Gotzmann's Persistence,
  and it is quite clear that we do not use the full strength of these two Theorems.
We believe that our proof can be generalised to cases when $d$ is much smaller with respect to $r$,
  but the number of cases to consider grows when we decrease $d$.
Thus a smarter, uniform treatment of all the cases is desired to get an effective bound on $d$.

\appendix
\section{Why cactus?}\label{section_why_cactus}

In this appendix we explain the name introduced in this article: \emph{the cactus variety}.

\begin{figure}[htb]
\centering
\begin{tabular}{cc}
\begin{minipage}[t]{0.42\textwidth}
\includegraphics[width=\textwidth]{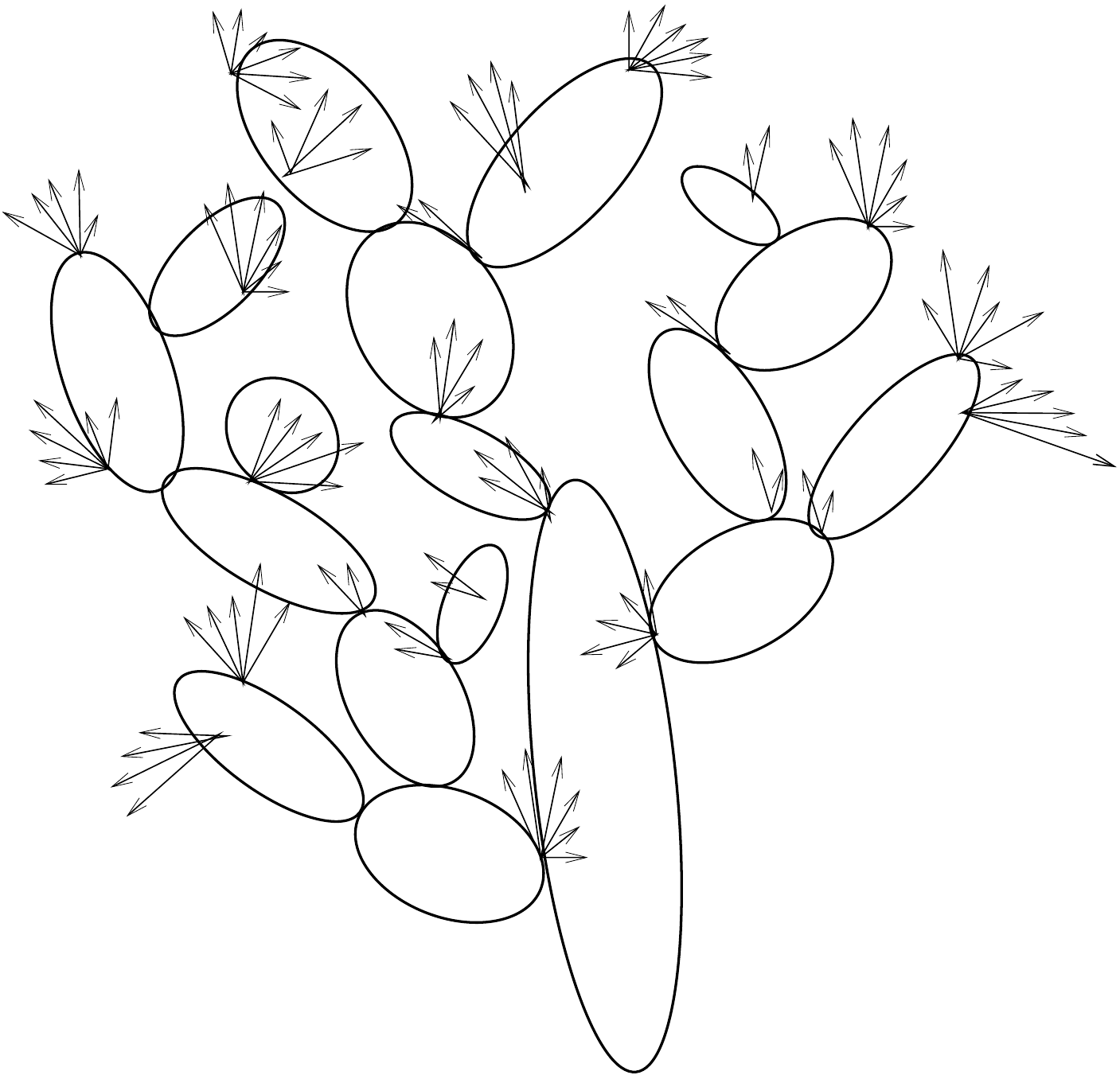}
\captionof{figure}{Our illustration of cactus variety.}
\label{figure_cactus_1}
\end{minipage}
&
\begin{minipage}[t]{0.52\textwidth}
\includegraphics[width=\textwidth]{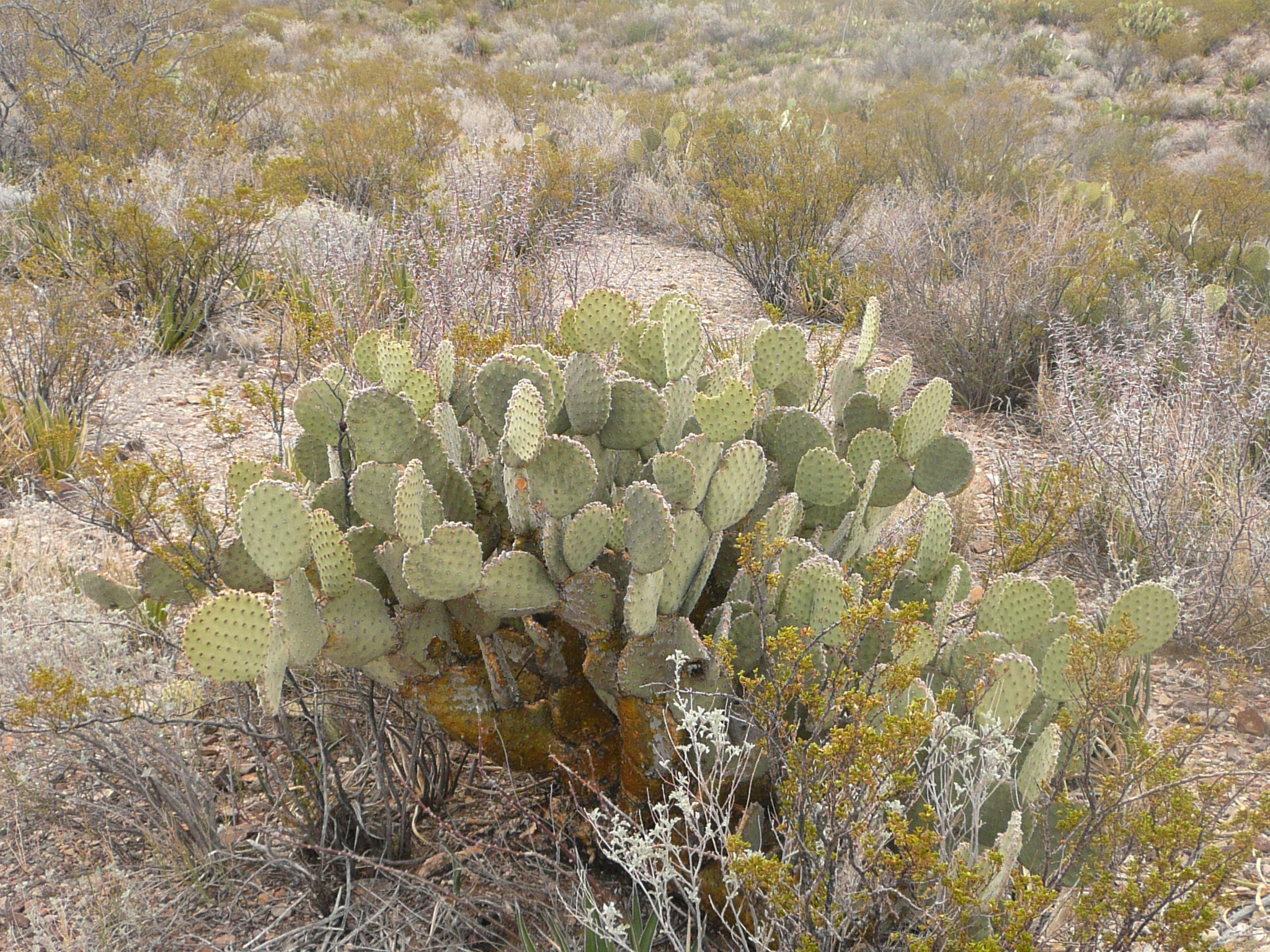}
\captionof{figure}{A picture of opuntia cactus in Big Bend National Park taken by the authors.}
\label{figure_cactus_2}
\end{minipage}
\end{tabular}
\end{figure}

The secant variety is swept by secant linear spaces, represented by ellipses on Figure~\ref{figure_cactus_1}.
These spaces form the \emph{stem} of the cactus.
The linear spans of non-smoothable Gorenstein schemes sometimes stick out of the secant variety,
  and these spans are the \emph{spines} of the cactus.

Our figure is very much simplified.
  We drew only a couple of tendentiously chosen secants, but there is an infinite number, a continuous family, of them.
  We drew $2$-dimensional secants (that is $r=3$ on the figure),
    whereas non-smoothable Gorenstein schemes exist only  for $r > 10$ or possibly even only for $r\ge 14$.
  The intersection of a secant plane and the linear span of a Gorenstein scheme might have positive dimension,
     rather than  dimension $0$, as might be suggested by the figure.
However, it is obvious that  some simplifications must be made,
   as it is impossible to draw a multidimensional object adequately in dimension 2.

Another illustrative comparison is the following:
   if one tries to wrap  a cactus stem in a tight and elegant package, the cactus spines might be an obstruction.
   We have obtained a similar obstruction (non-smoothable Gorenstein schemes)
     by trying to  present the defining equations of a secant variety
     in an easy and elegant form, that is as catalecticant minors.

We obtained Figure~\ref{figure_cactus_1} that reminds us of cacti of genus \emph{opuntia}
  --- see Figure~\ref{figure_cactus_2}.
These cacti the authors could observe in abundance in Texas,
  where they had recently spent two years.
Particularly, they could observe opuntia --- among many other exciting plants ---
  in Big Bend, a National Park in Wild West of Texas, during their Spring holiday.
On that trip, one cold and dark night,
  in a tent in Chisos Mountains, the second author started to realise that
  there might be a problem with the conjecture that all secant varieties to high degree Veronese reembeddings
  are cut out by catalecticant minors (see \cite[Question~1.2.2]{jabu_ginensky_landsberg_Eisenbuds_conjecture}).
After the holiday, the ideas were developed further in collaboration with Joseph Landsberg at Texas A\&M University.
Thus the name ``cactus variety'' is also our tribute to both: the inspiring wild life of Texas
  and the huge scientific and educational centre in Brazos County.

\bibliography{cactus_variety}
%\bibliography{references}
\bibliographystyle{alpha}

\end{document}